\documentclass[a4paper,11 pt]{article}

\usepackage{latexsym}
\usepackage[english]{babel}
\usepackage[utf8]{inputenc}  %<--per la tastiera italiana e le lettere accentate
\usepackage{amsmath} %<--pacchetto American Mathematician Society
\usepackage{amssymb} %<--pacchetto più ampio per simboli matematici
\usepackage{mathrsfs}
\usepackage{graphicx}
\usepackage{ esint }
\usepackage[usenames]{color}
\usepackage{enumerate}
\usepackage{color}
\usepackage{dsfont}
\usepackage{bbold}

\def\epsilon{\varepsilon}

\newtheorem{lemm}{Lemma}[section]

\newtheorem{teor}[lemm]{Theorem}

\newtheorem{deff}[lemm]{Definition}
\newtheorem{oss}[lemm]{Remark}

\def\dim{{\bf Proof. }}

\def\a1{$a_{1}$}
\def\r1{$r_{1}$}

\def\1{\{}
\def\2{\}}

\newcommand{\cvd}{\begin{flushright}$\Box$\end{flushright}}

\newcommand{\eq}{\begin{equation}}
\newcommand{\feq}{\end{equation}}
\newcommand{\be}{\begin{equation}}
\newcommand{\ee}{\end{equation}}

\begin{document}
\begin{center}

\vspace{1.cm}

{\Huge
\textbf{Stochastic Duality and \\ Orthogonal Polynomials}}

\vspace{0.5cm}

\begin{center}

{\Large
Chiara Franceschini$^{(1,2)}$, Cristian Giardin\`a$^{(1)}$}
\vspace{0.5cm}

\textit{$^{(1)}$ 
University of Modena and Reggio Emilia, via G. Campi 213/b, 41125 Modena, Italy.}

\textit{$^{(2)}$ 
University of Wisconsin-Madison, 480 Lincoln Drive, Madison (WI) 53706, USA.}

\end{center}

\end{center}

\vspace{1.cm}
\begin{center}
To Charles M. Newman on his $70^{th}$ birthday
\end{center}

\vspace{1.cm}

\begin{abstract}
\noindent
For a series of Markov processes we prove stochastic duality relations with
duality functions given by orthogonal polynomials. This means that expectations with 
respect to the original process (which evolves the variable of the orthogonal polynomial) 
can be studied via expectations with respect to the dual process (which evolves the index
of the  polynomial). The set of processes include interacting particle systems,
such as the exclusion process, the inclusion process and independent random walkers, as well as
interacting diffusions and redistribution models of Kipnis-Marchioro-Presutti type.
Duality functions are given in terms of classical orthogonal polynomials, both of discrete and 
continuous variable, and the measure in the orthogonality relation coincides with the process
stationary measure.
\end{abstract}

\newpage

\section{Motivations and main results}

\subsection{Introduction}

Duality theory is a powerful tool to deal with stochastic Markov processes
by which information on a given process can be extracted from another process,
its {\em dual}. The link between the two processes is provided by a set 
of so-called {\em duality functions}, i.e. a set of observables that are functions
of both processes and whose expectations, with respect to the two randomness,
can be placed in a precise relation (see Definition \ref{dop} below).
It is the aim of this paper to enlarge the space of duality functions
for a series of Markov processes that enjoy the stochastic duality
property. These novel duality functions are in turn connected to
the  polynomials that are orthogonal with respect to the inner product
provided by the stationary measure of those processes.
Thus the paper develops, via a series of examples, a connection between
probabilistic objects (duality functions) and  orthogonal polynomials.
See \cite{CED} for a previous example of a process with duality functions
in terms of Hermite polynomials.

Before going to the statement of the results, it is worth to stress the importance of a duality relation. Duality theory has been used in several different contexts.
Originally introduced for interacting particle systems in \cite{Spi70}
and further developed in \cite{Liggett}, the literature on stochastic duality
covers nowadays a host of examples. Duality has been applied -- among the others -- to boundary driven and/or bulk driven models of transport \cite{Spohn, BO, ImSa11, Ohku16}, Fourier's law
\cite{KMP, GKR},
diffusive particle systems and their hydrodynamic limit
\cite{Liggett,  DP06}, asymmetric interacting particle systems
scaling to KPZ equation \cite{GwaSpo92, BorCor13, BorCorSa14,  BCPS15, Cor15, BarCor16, CST16}, 
six vertex models \cite{BorCorGor16, CorPet16}, multispecies particle models
\cite{BelSch15,BelSch16, Kuan15, Kuan16, Kuan17}, correlation inequalities \cite{GRV10},
mathematical population genetics  \cite{M, CGGR}.
In all such different contexts it is used, in a way or another, the core simplification that duality provides.
Namely, in the presence of a dual process,  computation of $k$-point correlation functions of the
original process is mapped into the study of the evolution of only $k$ dual particles,
thus substantially reducing the difficulty of the problem. 

Besides applications, it is interesting to understand the mathematical structure behind duality.
This goes back to classical works by Sch\"utz and collaborators  \cite{Schutz-Sandow94, Sch97},
where the connection between stochastic duality and symmetries of quantum spin chains
was pointed out. More recently, the works \cite{GKRV,CGRT,CGRS16, Kuan16, Kuan17} 
further investigate this  framework and provide an algebraic approach to Markov processes with
duality  starting from a Lie algebra in the symmetric case, and its quantum deformation
in the asymmetric one. In this approach duality functions emerge as the intertwiners between
two different representations. Therefore one has a constructive theory, in which duality functions
arise from representation theory. 

An interesting problem is to fully characterize the set of all duality functions.
This question was first asked in \cite{M} where it was defined
the concept of duality space, i.e. the subspace of all  measurable functions 
on the configuration product space of two Markov processes for which
the duality relation (see Definition \ref{dop} below) holds. In \cite{M}
the dimension of this space is computed for some simple systems
and, as far as we know, a general answer is not available in the  
general case. Although the algebraic approach recalled above yields a duality function
from two representations of a (Lie) algebra, it is not clear a-priori if every duality 
function can be derived from this approach. 
In this paper we follow a different route. We shall show that duality functions can be placed in 
relation to orthogonal polynomials. 

\subsection{Results}

We recall the classical definition of stochastic duality.

\begin{deff}[Duality of processes.]
\label{dop}
Let $X=(X_{t})_{t\geq 0}$ and $N=(N_{t})_{t\geq 0}$ be two continuous time Markov processes with state spaces $\Omega$ and $\Omega^{dual}$ , respectively. We say that $N$ is {\em dual} to $X$ with duality function $D:\Omega\times \Omega^{dual}\longmapsto \mathbb{R} $ if 
\begin{equation}\label{du}
\mathbb{E}_{x}[D(X_{t},n)]=\mathbb{E}_{n}[D(x,N_{t})] \ ,
\end{equation}
for all $x \in \Omega$, $n \in \Omega^{dual}$ and $t\geq 0$.
In \eqref{du} $ \mathbb{E}_{x} $ (respectively $ \mathbb{E}_{n} $) is the expectation w.r.t. the law of the $X$ process initialized at $ x $ (respectively the $N$ process initialized at $n$).
If $X$ and $N$ are the same process, we say that $X$ is {\em self-dual} with self-duality function $D$.
\end{deff}
Under suitable hypothesis (see \cite{JK}), the above definition is equivalent to the definition of duality between Markov generators.
\begin{deff}[Duality of generators.] 
\label{dual2}
Let $L$ and $L^{dual} $ be generators of the two Markov processes $X=(X_{t})_{t\geq 0}$ and $N=(N_{t})_{t\geq 0}$, respectively.
We say that $L^{dual}$ is {\em dual} to $L$ with duality function $D:\Omega \times \Omega^{dual}\longrightarrow \mathbb{R}$ if 
\begin{equation}\label{duall}
[LD(\cdot,n)](x)=[L^{dual}D(x,\cdot)](n)
\end{equation}
where we assume that both sides are well defined.
In the case $L=L^{dual}$ we shall say that the process is {\em self-dual} and the self-duality relation becomes
\begin{equation}\label{self-duall}
[LD(\cdot,n)](x)=[LD(x,\cdot)](n)\;.
\end{equation}

\end{deff}
In \eqref{duall} (resp. \eqref{self-duall}) it is understood that $L$ on the lhs acts on $D$ as a function of the first variable $x$, while $L^{dual}$ (resp. $L$) on the rhs acts on $D$ as a function of the second variable $n$.
The definition \ref{dual2} is easier to work with, so we will always work under the assumption that the notion of duality (resp. self-duality) is the one in equation \eqref{duall} (resp. \eqref{self-duall}).

\begin{oss}
If $D(x,n)$ is duality function between two processes and $c:\Omega \times \Omega^{dual}\longrightarrow \mathbb{R}$ is constant under the dynamics of the two processes then $c(x,n)D(x,n)$ is also duality function. %In particular, examples of $c(x,n)$ could be the constant function or, for those processes that conserve the total number of particles, a function depending on this quantity.
We will always consider duality functions modulo the quantity $c(x,n)$.
For instance, the interacting particle systems studied in Section \ref{selfduality} conserve the total number of particles and thus $c$ is an arbitrary function of such conserved quantity.
\end{oss}

\vspace{0.3cm}

In this paper we shall prove that several Markov processes, for which a duality function is known from the algebraic approach, also admit 
a different duality function given in terms of polynomials that are orthogonal with respect to the process stationary measure. 

The processes that we consider include
discrete interacting particle systems (exclusion process, inclusion process and independent random walkers process)
as well as interacting diffusions (Brownian momentum process, Brownian energy process) and redistribution models that are obtained via a thermalization limit (Kipnis-Marchioro-Presutti processes). Their generators, that are defined in sections \ref{selfduality} and  \ref{duality},
have an algebraic structure from which duality functions have been previously derived \cite{GKRV}. 

The orthogonal polynomials we use are some of those with hypergeometric structure. 
More precisely we consider classical orthogonal polynomials, both discrete and continuous, 
with the exception of discrete Hahn polynomials and continuous Jacobi polynomials. Through this paper we follow the definitions of orthogonal polynomials
given in \cite{NSU}.

The added value of linking duality functions to orthogonal polynomials lies on the fact that they constitute an orthogonal basis of the associated Hilbert space.
Often in applications \cite{BorCorSa14, ImSa11, CGRT} some quantity of interest are expressed in terms of duality functions, for instance the current in interacting particle systems. This is then used in the study of the asymptotic properties and relevant scaling limits. For these reasons it seems reasonable that having an orthogonal basis of polynomials should be useful in those analysis.

\vspace{0.6cm}

The following theorem collects the results of this paper, details and rigorous proofs can be found in section \ref{selfduality} for self-duality and 
 section \ref{duality} for duality.

\begin{teor}
\label{teo}
For the processes listed below, the following duality relations hold true
 
%\noindent
%\begin{center}
%{\bf Self-duality}
%\bigskip
%
%\begin{tabular}{c|c|c}
%\it{Process} & \it{Stationary measure} & \it{Duality function}\\
%& & \it{in terms of}\\
%\hline 
%&&\\
%Exclusion Process & Binomial$(2j,p)$  &Krawtchouk polynomials\\
%with up to $2j$ particles &   &\\
%&&\\
%Inclusion Process & \qquad Negative Binomial $(2k,p)$ \qquad &Meixner polynomials\\
%with parameter $k$ &  &\\
%&&\\
%Independent  Walkers &Poisson($\lambda$) &Charlier polynomials\\
%\end{tabular}
%\end{center}
%
%\vspace{.5cm}
%
%\begin{center}
%
%{\bf Duality}
%\bigskip
%
%\begin{tabular}{c|c|c|c}
%\it{Process}&\it{Stationary measure}&\it{Duality function}&\it{Dual Process}\\
%&&\it{in terms of}&\\
%\hline
%&&&\\
%Brownian Momentum  & Gaussian$(0,\sigma^2)$ & Hermite polynomials & Inclusion process\\
%Process &  &  & with  $k=1/4$\\
%&&&\\
%Brownian Energy  & Gamma$(2k,\theta)$ &  Laguerre polynomials & Inclusion process\\
% Process with parameter $k$ &  &  & with parameter $k$\\
%&&&\\
%Kipnis-Marchioro-Presutti & Gamma$(2k,\theta)$ &  Laguerre polynomials & dual-KMP process\\
%with parameter $k$ &  &   & with parameter $k$\\
%\end{tabular}
%\end{center}
%

\begin{enumerate}

\item[i)] {\bf Self-duality}
\noindent
\begin{center}
\begin{tabular}{c|c|c}
\it{Process} & \it{Stationary measure} & \it{Duality function $D(x,n)$:}\\
& & \it{product of}\\
\hline 
&&\\
Exclusion Process & Binomial$(2j,p)$  &$K_n(x) / {2j \choose n}$\\
with up to $2j$ particles, SEP($j$) &  & \\
&&\\
Inclusion Process &  Negative Binomial $(2k,p)$  & $M_n(x) \frac{\Gamma(2k)}{\Gamma(2k+n)}$\\
with parameter $k$, SIP($k$) &  &\\
&&\\
Independent  Random &Poisson($\lambda$) &$C_n(x)$\\
Walkers, IRW & &\\
\end{tabular}

\vspace{0.3cm}
\noindent
where $K_n(x)$ stands for Krawtchouk polynomials, $M_n(x)$ for Meixner polynomials,  $C_n(x)$ for Charlier polynomials.
\end{center}

\vspace{.5cm}

\item[ii)] {\bf Duality}
\noindent
\begin{center}
\begin{tabular}{c|c|c|c}
\it{Process}&\it{Stationary measure}&\it{Duality function} $D(x,n)$&\it{Dual Process}\\
&&\it{product  of}&\\
\hline
&&&\\
Brownian Momentum  & Gaussian$(0,\sigma^2)$ & $\frac{1}{(2n-1)!!} H_{2n}(x)$ & Inclusion process\\
Process, BMP &  &  & with  $k=1/4$\\
&&&\\
Brownian Energy Process & Gamma$(2k,\theta)$ &   $  \frac{n! \Gamma(2k)}{\Gamma(2k+n)}L_n^{(2k-1)}(x)$ & Inclusion process\\
with parameter $k$, BEP($k$)  &  &  & with parameter $k$\\
&&&\\
 Kipnis-Marchioro-Presutti  & Gamma$(2k,\theta)$ &  $  \frac{n! \Gamma(2k)}{\Gamma(2k+n)} L_n^{(2k-1)}(x) $ & dual-KMP($k$) process\\
with parameter $k$,  KMP($k$) &  &   & with parameter $k$\\
\end{tabular}

\vspace{0.3cm}
\noindent
where $H_n(x)$ stands for Hermite polynomials and $L^{(2k-1)}_n(x)$ for generalized Laguerre polynomials.
\end{center}

\end{enumerate}

\end{teor}

%\begin{center}
%\begin{tabular}{rcc}
%\textbf{Process}&\textbf{Stationary measure}&\textbf{Duality function in terms of}\\
%SEP($j$)&Bin($2j$,$p$)&Krawtchouk polynomials\\
%SIP($k$)&NegBin($2k$,$p$)&Meixner polynomials\\
%IRW&Poi($\lambda$)&Charlier polynomials\\
%\end{tabular}
%\end{center}
%
%
%
%\begin{center}
%\begin{tabular}{rccc}
%\textbf{Process}&\textbf{Stationary measure}&\textbf{Duality function in terms of}&\textbf{Dual Process}\\
%BMP&N($0$,$1/2$)&Hermite polynomials& SIP($1/4$)\\
%BEP($k$)&Gamma($2k$,$\theta$)&Generalized Laguerre polynomials& SIP($k$)\\
%KMP($k$)&Gamma($2k$,$\theta$)&Generalized Laguerre polynomials& SIP($k$)\\
%\end{tabular}
%\end{center}

\vspace{.3cm}
\begin{oss}
As can be inferred from the table, the duality function is not, in general, the orthogonal polynomial itself, but a suitable normalization. 
Moreover, for the process defined on multiple sites the duality functions exhibit a product structure 
where each factor is in terms of the relevant orthogonal polynomial.
\end{oss}

\subsection{Comments}

\paragraph{Old and new dualities.}
It is known that the six processes we consider satisfy the same (self-)duality relation described by the chart above with different (self-)duality functions.
New and old (self-)duality functions can be related using the explicit form of the polynomial that appears in the new one via a relation of the
following type:
\begin{equation}
D^{new}(x,n)=\sum_{k=0}^{n}d(k,n)D^{old}(x,k)
\end{equation} 
where $d(k,n)$ also depends on the parameter of the stationary measure.
If, for example, we consider the self-duality of the Independent Random Walker, where the
new self-duality functions are given by the Charlier polynomials and the old self-duality functions are given by
the falling factorial, then we have
\be
D^{new}(x,n) = \sum_{k=0}^{n} \binom {n} {k} (-\lambda)^{n-k} \dfrac{x!}{(x-k)!}   \qquad\qquad D^{old}(x,n)=\dfrac{x!}{(x-n)!}
\ee
so that
\begin{equation}
D^{new}(x,n)=\sum_{k=0}^{n} \binom {n} {k} (-\lambda)^{n-k} D^{old}(x,k)
\end{equation} 
where $\lambda$ is the Poisson distribution parameter.

\paragraph{The norm choice.}
A sequence of orthogonal polynomials   $\{p_n(x) : n\in \mathbb{N}\}$ on the interval $(a,b)$ is defined by the 
choice of a measure $\mu$ to be used in the scalar product 
and a choice of a norm $d_n^2>0$
\be
\langle p_n, p_m\rangle := \int_{a}^b p_n(x) p_m(x) d\mu(x) = \delta_{n,m} d_n^2\;.
\ee
The main result of this paper states that orthogonal polynomials with a properly chosen {\em normalization} 
furnish duality functions
for the processes we consider. 
In our setting the probability measure $\mu$ to be used is provided by the marginal 
stationary measure of the process. The appropriate norm has been found by making the ansatz
that new duality functions can be obtained by the Gram-Schmidt orthogonalization 
procedure  initialized with the old duality functions derived from the algebraic approach \cite{GKRV}.
Namely, if one starts from a sequence $p_n(x) $ of orthogonal polynomials with norm $d_n^2$
and claims that $D^{new}(x,n) = b_n p_n(x)$, where $b_n$ is the appropriate normalization, then the previous ansatz yields
\be\label{normalizz}
b_n = \frac{\langle D^{old}(\cdot,n), p_n\rangle}{d_n^2}\;.
\ee

\paragraph{General strategy of the proof.}
Once a properly normalized orthogonal polynomial is identified as a candidate duality function, then 
the proof of the duality statement is obtained via explicit computations.
Each proof heavily relies on the {\em hypergeometric structure} of the polynomials involved. 
The idea is to use three structural properties of the polynomials family $p_n(x) $, i.e. the differential or difference hypergeometric equation they satisfy, their three terms recurrence relation and 
the expression for the raising ladder operator. Those properties are then transported to the duality function using the proper 
rescaling  $b_n p_n(x)$. This yields three identities for the duality function, that can be used in the expression of the process generator. Finally, algebraic manipulation allows to verify relation \eqref{duall} for duality and \eqref{self-duall} for self-duality.

\paragraph{Stochastic duality and polynomials duality}
Last, we point out that there exists a notion of duality within the context of orthogonal polynomials, see Definition 3.1 in  \cite{KLS}.
For example, in case of discrete orthogonal polynomials,  the polynomials self-duality can be described by the identity $p_{n}(x)=p_{x}(n)$ where $x$ and $n$
take values in the same discrete set. It turns out that Charlier polynomials are self-dual, whereas Krawtchouk and Meixner polynomials, defined as in \cite{NSU} see also section \ref{selfduality}, are not. However, it is possible to rescale them by a constant $b_{n}$ so that they become self-dual. The constant $b_{n}$ is indeed provided by \eqref{normalizz}.
It is not clear if the two notions of stochastic duality and polynomials duality are in some relation one to the other.
This will be investigated in a future work.

\subsection{Paper organization }
The paper is organized as follows. Section \ref{hpoly} is devoted to a (non exhaustive) review of hypergeometric orthogonal polynomials, 
we recall their key properties that are crucial in proving our results. We closely follow \cite{NSU} and the expert reader might skip this part
without being affected. 
The original results are presented in sections \ref{selfduality} and \ref{duality}. In section \ref{selfduality} we describe three interacting particle 
systems that are self-dual and prove the statement of Theorem \ref{teo}, part  i).
Section \ref{duality} is dedicated to the duality relations of two diffusion processes and a jump process obtained as thermalization limit 
of the previous ones and it contains the proof of Theorem \ref{teo}, part ii).

\section{Preliminaries: hypergeometric orthogonal polynomials}\label{hpoly}
In this section we give a quick overview of the continuous and the discrete hypergeometric polynomials (see \cite{TSC, KLS, NSU, W}) by reviewing some of their structural properties that will be used in the following. 

We start by recalling that the hypergeometric orthogonal polynomials arise from an hypergeometric equation, whose solution can be written in terms of an hypergeometric function $\mathstrut_r F_s$.
\begin{deff}[Hypergeometric function]
The hypergeometric function is defined by the series
\begin{equation}\label{hyf}
\mathstrut_r F_s \left( {\left. \genfrac{}{}{0pt}{} {a_{1},\ldots,a_{r} } { b_{1},\ldots,b_{s} }  \right\vert {x}} \right) = \sum_{k=0}^{\infty}\dfrac{(a_{1})_{k}\cdots(a_{r})_{k}}{(b_{1})_{k}\cdots(b_{s})_{k}}\dfrac{x^{k}}{k!}
\end{equation}
where $( a )_{k}$ denotes the Pochhammer symbol defined in terms of the Gamma function as
\begin{equation*}
\left( a \right)_{k}=\dfrac{\Gamma(a+k)}{\Gamma(a)}.
\end{equation*}
\end{deff}
\begin{oss}
\label{poly-geo}
Whenever one of the numerator parameter $a_{j}$ is a negative integer $-n$, the hypergeometric function  $\mathstrut_r F_s $ is a finite sum up to $n$, i.e. a polynomial in $x$ of degree $n$.
\end{oss}

\noindent
{\large\textbf{The continuous case.}} Consider the hypergeometric differential equation 
\begin{equation}\label{hdf}
\sigma(x)y^{''}(x)+\tau(x)y^{'}(x)+\lambda y(x)=0
\end{equation}
where $\sigma(x)$ and $\tau(x)$ are polynomials of at most second and first degree respectively and $\lambda$ is a constant. 
%Its solutions are called functions of hypergeometric type. 
A peculiarity of the hypergeometric equation is that, for all $n$, $y^{(n)}(x)$, i.e. the $n^{th}$ derivative of a solution $y(x)$, also  solves 
an hypergeometric equation, namely 
\begin{equation}\label{hdf1}
\sigma(x)y^{(n+2)}(x)+\tau_{n}(x)y^{(n+1)}(x)+\mu_{n} y^{(n)}(x)=0
\end{equation}
with 
\begin{equation}\label{deftau}
\tau_{n}(x)=\tau(x)+n\sigma^{'}(x) 
\end{equation} 
and
\begin{equation}
\mu_{n}=\lambda +n\tau^{'}+\frac{1}{2}n(n-1)\sigma^{''}\;.
\end{equation}
%is a constant.\\
%To show this, let's start by differentiating equation \eqref{hdf} once 
%\begin{equation*}
%\sigma(x)^{'}y(x)^{''}+\sigma(x)y(x)^{'''}+\tau(x)^{'}y(x)^{'}+\tau(x)y(x)^{''}+\lambda  y(x)^{'}=0.
%\end{equation*}
%Rearranging leads to
%\begin{equation*}
%\sigma(x)y(x)^{'''}+(\sigma(x)^{'}+\tau(x))y(x)^{''}+ 
%(\tau(x)^{'}+\lambda ) y(x)^{'}=0
%\end{equation*}
%and calling $\tau_{1}=\sigma(x)^{'}+\tau(x)$ and $\mu_{1}=\tau(x)^{'}+\lambda$ let us recover the hypergeometric structure for a differential equation whose solution is $y^{'}$. $n$ times iteration leads to equation \eqref{hdf1}.\\
We concentrate on a specific family of solutions: for each $n \in \mathbb{N}$, let $\mu_{n}=0$, so that 
\begin{equation} \label{defln1}
\lambda =\lambda_{n}=-n\tau^{'}-\frac{1}{2}n(n-1)\sigma^{''}
\end{equation}
and equation \eqref{hdf1} has a particular  solution given by $y^{(n)}(x)$ constant. This implies that $y(x)$ is a polynomial of degree $n$, called \emph{polynomial of hypergeometric type} 
(see remark \ref{poly-geo}) and denoted by $p_{n}(x)$.
In the following we will assume that those polynomials are of the form
\be
\label{leading}
p_{n}(x)=a_{n}x^{n}+b_{n}x^{n-1}+\ldots \qquad \qquad \qquad a_{n}\neq 0\;.
\ee
It is well known \cite{NSU} that polynomials of hypergeometric type satisfy the orthogonality relation
\begin{equation}\label{or}
\int_{a}^{b}p_{n}(x)p_{m}(x) \rho(x) dx= \delta_{n,m}d^{2}_{n}(x)
\end{equation}
for some (possibly infinite) constants $a$ and $b$ and where the function $\rho(x)$ satisfies the differential equation
\begin{equation}\label{rhho}
(\sigma \rho)^{'}=\tau \rho\;.
\end{equation}
The sequence $d^{2}_{n}$ can be written in terms of $\sigma(x), \rho(x)$ and $a_n$ as
\begin{equation}
d_{n}^{2}=\dfrac{(a_{n} n!)^{2}}{\prod_{k=0}^{n-1}(\lambda_n - \lambda_k)} \int_{a}^{b}(\sigma(x))^n\rho(x) dx.
\end{equation}
As a consequence of the orthogonal property the polynomials of hypergeometric type satisfy a three terms recurrence relation
\begin{equation}\label{expan2-cont}
xp_{n}(x)=\alpha_{n}p_{n+1}(x)+\beta_{n}p_{n}(x)+\gamma_{n}p_{n-1}(x)
\end{equation}
where 
\be
\alpha_{n}=c_{n+1,n} \qquad \beta_{n}=c_{n,n} \qquad \gamma_{n}=c_{n-1,n}
\ee
with
\be
c_{k,n} = \dfrac{1}{d_{k}^{2}}\int_{a}^{b}p_{k}(x)x p_{n}(x) \rho(x) dx \;.
\ee
The coefficients $\alpha_n,\beta_n,\gamma_n$ can be expressed in terms of the squared norm $d_n^2$ and the leading coefficients $a_n, b_n$ in \eqref{leading} as \cite{NSU}
\begin{equation}
\alpha_{n}=\dfrac{a_{n}}{a_{n+1}} \qquad \beta_{n}=\dfrac{b_{n}}{a_{n}}-\dfrac{b_{n+1}}{a_{n+1}} \qquad\qquad \gamma_{n}=\dfrac{a_{n-1}}{a_{n}}\dfrac{d_{n}^{2}}{d_{n-1}^{2}}\;.
\end{equation}
Finally, we will use the raising operator $\mathrm{R}$ that, acting on the polynomials $p_{n}(x)$, provides the polynomials of degree $n+1$.
Such an operator is obtained from the Rodriguez formula, which provides an explicit form for polynomials of hypergeometric type
\begin{equation}\label{rod}
p_{n}(x)=\dfrac{B_{n}( \sigma^{n}(x) \rho(x))^{(n)}}{\rho(x)} \qquad\qquad \text{with} \qquad\qquad B_{n}=\dfrac{a_n}{ \prod_{k=0}^{n-1}\left(  \tau^{'}+\frac{n+k-1}{2} \sigma^{''}\right)} \;.
\end{equation}
The  expression of  the raising operator  (see eq. 1.2.13 in \cite{NSU}) reads
\be
\label{creation}
\mathrm{R} p_n(x)  = r_n p_{n+1}(x)
\ee
where
\begin{equation}
\mathrm{R}p_{n}(x)=  \lambda_{n}\tau_{n}(x) p_n(x) - n \sigma(x) \tau_{n}^{'} p_{n}^{'}(x) \qquad \qquad \text{and} \qquad r_n= \lambda_{n}\dfrac{B_{n}}{B_{n+1}}\;.
\end{equation}
Remark that the raising operator increases the degree of the polynomial by one, similarly to the so-called  backward shift operator \cite{KLS}. However the raising operator in \eqref{creation} does not change the parameters involved in the function $\rho$, whereas the backward operator increases the degree and lowers the parameters \cite{Koornwinder}.

\vspace{.5cm}
\noindent{\large\textbf{The discrete case.}} Everything discussed for the continuous case has a discrete analog, where the derivatives are replaced by the discrete difference derivatives. In particular it is worth mentioning that
\begin{equation*}
\Delta f(x)=f(x+1)-f(x) \qquad \text{and}  \qquad \nabla f(x)=f(x)-f(x-1)\;.
\end{equation*} 
The corresponding hypergeometric differential equation \eqref{hdf} is the discrete hypergeometric difference equation
\begin{equation}\label{hyperdiscrete}
\sigma(x)\Delta \nabla y(x) + \tau(x) \Delta y(x) + \lambda y(x)=0
\end{equation}
where $\sigma(x)$ and $\tau(x)$ are polynomials of second and first degree respectively, $\lambda$ is a constant. The differential equation solved by the $n^{th}$ discrete derivative of $y(x)$, $y^{(n)}(x):=\Delta^{n}y(x)$, is the solution of another difference equation of hypergeometric type
\begin{equation}\label{yn}
\sigma(x)\Delta \nabla y^{(n)}(x) + \tau_{n} \Delta y^{(n)}(x) + \mu_{n} y^{(n)}(x)=0
\end{equation}
with
\begin{equation}
\tau_{n}(x)=\tau(x+n)+\sigma(x+n)- \sigma(x)
\end{equation}
and
\begin{equation}
\mu_{n}=\lambda+n\tau^{'}+\frac{1}{2}n(n-1)\sigma^{''}\;.
\end{equation}
If we impose $\mu_{n}=0$, then
\begin{equation} \label{defln2}
\lambda =\lambda_{n}=-n\tau^{'}-\frac{1}{2}n(n-1)\sigma^{''}
\end{equation}
and $y^{(n)}(x)$ is a constant solution of equation \eqref{yn}. Under these conditions, $y(x)$, solution of \eqref{hyperdiscrete},  is a polynomial of degree $n$, called \emph{discrete polynomial of hypergeometric type} 
(see remark \ref{poly-geo}) and denoted by $p_{n}(x)$.\\
%In the following we will assume that those polynomials are of the form
%\be
%\label{leading}
%y_{n}(x)=a_{n}x^{n}+b_{n}x^{n-1}+\ldots
%\ee
The derivation of the orthogonal property is done in a similar way than the one for the continuous case where the integral is replaced by a sum
\begin{equation}
\sum_{x=a}^{b-1}p_{n}(x)p_{m}(x)\rho(x)=\delta_{n,m}d_{n}^{2}
\end{equation}
constants $a$ and $b$ can be either finite or infinite and the function $\rho(x)$ is solution of
\begin{equation}
\Delta[\sigma(x)\rho(x)]=\tau(x)\rho(x)\;.
\end{equation}
The sequence $d_{n}^{2}$ can be written in terms of $\sigma(x), \rho(x)$ and $a_{n}$ as
\begin{equation}
d_{n}^{2}=\dfrac{(a_{n}n!)^{2}}{\prod_{k=0}^{n-1}(\lambda_{n}-\lambda_{k})} \sum_{x=a}^{b-n-1} \left( \rho(x+n)\prod_{k=1}^{n} \sigma(x+k) \right) 
\end{equation}
As a consequence of the orthogonal property, the discrete polynomials of hypergeometric type satisfy a three terms recurrence relation
\begin{equation}\label{expan2}
xp_{n}(x)=\alpha_{n}p_{n+1}(x)+\beta_{n}p_{n}(x)+\gamma_{n}p_{n-1}(x)
\end{equation}
where 
\be
\alpha_{n}=c_{n+1,n} \qquad \beta_{n}=c_{n,n} \qquad \gamma_{n}=c_{n-1,n}
\ee
with
\be
c_{k,n} = \dfrac{1}{d_{k}^{2}}\sum_{x=a}^{b}p_{k}(x)x p_{n}(x)\;.
\ee
The coefficients $\alpha_n,\beta_n,\gamma_n$ can be expressed in terms of the squared norm $d_n^2$ and the leading coefficients $a_n, b_n$ in \eqref{leading} as \cite{NSU}
\begin{equation}
\alpha_{n}=\dfrac{a_{n}}{a_{n+1}} \qquad \beta_{n}=\dfrac{b_{n}}{a_{n}}-\dfrac{b_{n+1}}{a_{n+1}} \qquad\qquad \gamma_{n}=\dfrac{a_{n-1}}{a_{n}}\dfrac{d_{n}^{2}}{d_{n-1}^{2}}\;.
\end{equation}
The discrete Rodriguez formula
\begin{equation}
p_{n}(x)=\dfrac{B_{n}}{\rho(x)} \nabla^{n}\left[ \rho(x+n) \prod_{k=1}^{n} \sigma(x+k) \right]  \qquad\qquad \text{with} \qquad\qquad  B_{n}=\dfrac{a_{n} }{ \prod_{k=0}^{n-1}\left( \tau^{'}+\frac{n+k-1}{2}\sigma^{''} \right) }
\end{equation}
leads to an expression for the discrete raising operator $\mathrm{R}$ (see eq. 2.2.10 in \cite{NSU})
\begin{equation}
\label{raising-discrete}
\mathrm{R}p_{n}(x)=r_{n}p_{n+1}(x)
\end{equation}
where 
\begin{equation}
\mathrm{R}p_{n}(x)= \left[ \lambda_{n} \tau_{n}(x) -n\tau_{n}^{'} \sigma(x) \nabla \right]p_{n}(x) \qquad \qquad r_{n}=\lambda_{n}\dfrac{B_{n}}{B_{n+1}} \;.
\end{equation}
We remark again that the raising operator shouldn't be confused with the backward shift operator in \cite{Koornwinder} which changes the value of parameters of the distribution $\rho$.

\section{Self-duality: proof of Theorem \ref{teo} part i)}\label{selfduality}
We consider first self-duality relations for some discrete interacting particle systems. 
In this case the dual process is an independent copy of the original one.                                                                                                                                                                                                            
Notice that, even if the original process and its dual are the same, a massive simplification occurs, namely
the $n$-point correlation function of the original process can be expressed by duality in terms of
only $n$ dual particles. Thus a problem for many particles, possibly infinitely many in the infinite volume,
may be studied via a finite number of dual walkers.

\subsection{The  Symmetric Exclusion Process, SEP($j$)}
In the  Symmetric Exclusion Process with parameter $j\in\mathbb{N}/{2}$, denoted by SEP($j$), each site can have at most $2j$ particles, 
jumps occur at rate proportional to the number of particles in the departure site times the number of holes in the arrival site.
The special case $j=1/2$ corresponds to the standard exclusion process with hard core exclusion, i.e. each site can be either full or empty \cite{Liggett}. We consider the setting where the spacial structure is given by the undirected connected graph $G=(V,E)$ with $N$ vertices and edge set $E$. A particle configuration is denoted by $\textbf{x}=(x_{i})_{i \in V}$ where $x_{i} \in {\{0,\ldots,2j\}}$ is interpreted as the particle number at vertex $i$. The process generator reads 
\begin{equation}\label{sepj}
L^{SEP(j)}f({\bf x})=\sum_{\substack{1\le i < l \le N   \\  (i,l)\in E}} x_{i}(2j-x_{l})\left[ f({\bf x}^{i,l})-f({\bf x})\right]+(2j-x_{i})x_{l}\left[ f({\bf x}^{l,i})-f({\bf x})\right]
\end{equation}
where ${\bf x}^{i,l}$ denotes the particle configuration obtained from the configuration ${\bf x}$ by moving one particle from vertex $i$ to vertex $l$
and
where $f: {\{0,1,\ldots,2j\}}^N \to \mathbb{R}$ is a function in the domain of the generator.

It is easy to verify that the reversible (and thus stationary) measure of the process for every $p\in(0,1)$ is given by the homogeneous product measure with marginals the Binomial distribution
with parameters $2j>0$ and  $p\in(0,1)$, i.e. with probability mass function
\be
\label{rho-bin}
\rho(x) = {2j \choose x} p^x (1-p)^{2j-x}\;, \qquad\qquad x\in\{0,1,\ldots,2j\}\;.
\ee
The orthogonal polynomials with respect to the Binomial distribution are the Krawtchouk polynomials $K_n(x)$ with parameter $2j$ \cite{Kra}. 
These polynomials are obtained by choosing in Eq. \eqref{hyperdiscrete}
\be
\sigma(x) = x \qquad\qquad \tau(x) = \frac{2j p -x}{1-p} \qquad\qquad \lambda_n = \frac{n}{1-p}\;.
\ee 
Equivalently, Krawtchouk polynomials are polynomial solutions of the finite difference equation
\be
\label{diff-eq-kra}
x [K_n(x+1) - 2 K_n(x) + K_n(x-1)] + \frac{2j p -x}{1-p} [K_n(x+1) - K_n(x)] + \frac{n}{1-p} K_n(x) = 0 \;. 
\ee
They satisfy the orthogonality relation
\be
\sum_{x=0}^{2j} K_n(x) K_m(x) \rho(x)= \delta_{n,m} d_n^2
\ee
with $\rho$ as in \eqref{rho-bin} and norm in $\ell^2(\{0,1,\ldots,2j\}, \rho)$ given by
\be
d_n^2 = {2j \choose n} p^n (1-p)^n\;.
\ee
As a consequence of the orthogonality they satisfy the recurrence relation \eqref{expan2}
with 
\be
\alpha_n = n+1 \qquad \qquad \beta_n = n+2jp-2np  \qquad \qquad \gamma_n = p (1-p)(2j-n+1)
\ee
and the three-point recurrence then becomes
\be
\label{recu-kra}
xK_n(x) = (n+1)K_{n+1}(x) + (n+2jp -2np)K_n(x) + p (1-p)(2j-n+1)K_{n-1}(x)\;.
\ee
Furthermore, the raising operator in \eqref{raising-discrete} provides the relation 
%\be
%a f_n(x)  = \frac{p}{1-p} (n+x-2j) f_n(x) + xf_n(x-1) 
%\ee
%with 
%\be
%c(n) =  \frac{n-1}{1-p}
%\ee
%and thus we have the equation
\be
\label{raising-kra}
xK_n(x-1)  + \frac{p}{1-p} (n+x-2j) K_n(x)  = \frac{n-1}{1-p}K_{n+1}(x)\;. 
\ee

\subsubsection{Self-duality for SEP($j$)}
The Krawtchouk polynomials $K_n(x)$  do not satisfy a self-duality relation in the sense of Definition \ref{dual2}. 
However, the following theorem shows that, with a proper normalization, it is possible to find a duality function related to them.
%The proof of the theorem  relies on the structural properties of the Krawtchouk polynomials.
\begin{teor}
The SEP($j$) is a self-dual Markov process with self-duality function 
\begin{equation}\label{kra}
D_{{\bf n}}({\bf x})= \prod_{i=1}^N \dfrac{n_i!(2j-n_i)!}{2j!}K_{n_i}(x_i)
\end{equation}
where $K_{n}(x)$ denotes the Krawtchouk polynomial of degree $n$.
\end{teor}
\begin{oss}
Since the Krawtchouk polynomials can be rewritten in terms of hypergeometric functions, the self-duality function turns out to be
\begin{equation}
D_{\bf n}({\bf x})= \prod_{i=1}^N \mathstrut_2 F_1 \left( {\left. \genfrac{}{}{0pt}{} {-n_i,-x_i } { -2j }  \right\vert {\frac{1}{p}}} \right) 
\end{equation}
where $\mathstrut_2 F_0$ is the hypergeometric function defined in \eqref{hyf}.
\end{oss}
\dim
We need to verify the self-duality relation in equation \eqref{self-duall}. Since the generator of the process is a sum of terms acting on two variables only, we shall verify the self-duality relation for two sites, say $1$ and $2$. We start by writing the action of the SEP($j$) generator working on the duality function for these two sites: 
\begin{eqnarray}
L^{SEP(j)}D_{n_{1}}(x_{1})D_{n_{2}}(x_{2}) 
& =  & x_{1}(2j-x_{2}) \left[ D_{n_{1}}(x_{1}-1)D_{n_{2}}(x_{2}+1) - D_{n_{1}}(x_{1})D_{n_{2}}(x_{2})\right]   \\
& +  & (2j-x_{1})x_{2} \left[ D_{n_{1}}(x_{1}+1)D_{n_{2}}(x_{2}-1) - D_{n_{1}}(x_{1})D_{n_{2}}(x_{2}) \right] \nonumber
\end{eqnarray}
rewriting this by factorizing site 1 and 2, i.e.
\begin{eqnarray}
\label{step2}
L^{SEP(j)}D_{n_{1}}(x_{1})D_{n_{2}}(x_{2}) 
& =  & x_{1}D_{n_{1}}(x_{1}-1)(2j-x_{2})D_{n_{2}}(x_{2}+1)-x_{1}D_{n_{1}}(x_{1})(2j-x_{2})D_{n_{2}}(x_{2})   \\ 
& + &(2j-x_{1})D_{n_{1}}(x_{1}+1) x_{2}D_{n_{2}}(x_{2}-1)- (2j-x_{1}) D_{n_{1}}(x_{1})x_{2}D_{n_{2}}(x_{2}) \nonumber
\end{eqnarray} 
we see that we need an expression for the following terms:
\be
\label{3-termini}
xD_{n}(x)\;, \qquad\quad xD_{n}(x-1)\;, \qquad\quad (2j- x)D_{n}(x+1)\;.
\ee
To get those we first write the difference equation \eqref{diff-eq-kra}, the recurrence relation \eqref{recu-kra} and the raising operator equation \eqref{raising-kra} in terms of $D_{n}(x)$. 
This is possible using equation \eqref{kra} that provides $D_n(x)$ as a suitable normalization of $K_n(x)$, i.e. 
\be
D_{n}(x)= \dfrac{n!(2j-n)!}{2j!}K_{n}(x)\;.
\ee   
Then the first term in \eqref{3-termini} is simply obtained from the normalized recurrence relation, whereas the second and third terms are provided by simple algebraic manipulation of the normalized raising operator equation and the normalized difference equation.
We get 
\begin{eqnarray}
xD_{n}(x)=-p(2j-n)D_{n+1}(x)+(n+2pj-2pn)D_{n}(x)-n(1-p)D_{n-1}(x)\\
xD_{n}(x-1)=-p(2j-n)D_{n+1}(x) +p(2j-2n)D_{n}(x)+npD_{n-1}(x)\\
(2j- x)D_{n}(x+1)=p(2j-n)D_{n+1}(x)+(1-p)(2j-2n)D_{n}(x)-\dfrac{n}{p}(1-p)^{2}D_{n-1}(x)\;. 
\end{eqnarray}
These expressions can now be inserted into \eqref{step2}, which then reads:
\begin{align*}
L^{SEP(j)}&D_{n_{1}}(x_{1})D_{n_{2}}(x_{2}) 
\\ &=
\left[ p(n_{1}-2j)D_{n_{1}+1}(x_1) + p(2j-2n_{1})D_{n_{1}}(x_1) + pn_{1}D_{n_{1}-1}(x_1)\right] \times 
\\ &\left[ p(2j-n_{2})D_{n_{2}+1}(x_{2}) +(1-p)(2j-2n_{2})D_{n_{2}}(x_{2}) -\dfrac{n_{2}}{p}(1-p)^{2}D_{n_{2}-1}(x_{2})  \right]  \\&
+\left[p(2j-{n_{1}})D_{n_{1}+1}(x_1) -(n_{1}+2jp-2pn_{1})D_{n_{1}}(x_1) +(1-p)n_{1}D_{n_{1}-1}(x_1) \right] \times \\ & 
\left[ p(2j-n_{2})D_{n_{2}+1}(x_2)-(n_{2}+2pj-2pn_{2})D_{n_{2}}(x_2)+n_{2}(1-p)D_{n_{2}-1}(x_2) +2jD_{n_{2}}(x_2)\right] \\& 
+\left[ p(n_{2}-2j)D_{n_{2}+1}(x_{2})  +p(2j-2n_{2})D_{n_{2}}(x_{2})  + pn_{2}D_{n_{2}-1}(x_{2})    \right]  \times \\ &
\left[ p(2j-n_{1})D_{n_{1}+1}(x_1)+(1-p)(2j-2n_{1})D_{n_{1}}(x_1)-\dfrac{n_{1}}{p}(1-p)^{2}D_{n_{1}-1}(x_1) \right]  \\ &
+\left[p(2j-{n_{2}})D_{n_{2}+1}(x_{2})  -(n_{2}+2pj-2pn_{2})D_{n_{2}}(x_{2})  + (1-p)n_{2}D_{n_{2}-1}(x_{2}) \right] \times \\& 
\left[ p(2j-n_{1})D_{n_{1}+1}(x_1)-(n_{1}+2pj-2pn_{1})D_{n_{1}}(x_1)+n_{1}(1-p)D_{n_{1}-1}(x_1) +2jD_{n_{1}}(x_1)\right]\;.
\end{align*}
Working out the algebra, substantial simplifications are revealed in the above expression.
A long but straightforward computation shows that only products of polynomials with degree $n_1+n_2$ survive.
%At this point it is suffice to notice that the coefficients of $D_{n_{1}+1}D_{n_{2}+1}$, $D_{n_{1}-1}D_{n_{2}-1}$, $D_{n_{1}+1}D_{n_{2}}$, $D_{n_{1}}D_{n_{2}+1}$, $D_{n_{1}-1}D_{n_{2}}$, $D_{n_{1}}D_{n_{2}-1}$, $D_{n_{1}-1}D_{n_{2}+1}$ and $D_{n_{1}+1}D_{n_{2}-1}$ are all zero. 
In particular, after simplifications, one is left with
%\begin{eqnarray*}
%L^{SEP(j)}D_{n_{1}}(x_1)D_{n_{2}}(x_2)
%& = &(2j-n_{1})D_{n_{1}+1}(x_1)n_{2}D_{n_{2}-1}(x_2)\left[p\dfrac{1}{p}(1-p)^{2}+p(1-p)+p^{2}+p(1-p) \right] \\ 
%& + & n_{1}D_{n_{1}-1}(x_1)(2j-n_{2})D_{n_{2}+1}(x_2)\left[p\dfrac{1}{p}(1-p)^{2}+p(1-p)+p^{2}+p(1-p) \right] \\ 
%& + & D_{n_{1}}(x_1)D_{n_{2}}(x_2)\left[ 2n_{1}n_{2}-2jn_{1}-2jn_{2}\right] 
%\end{eqnarray*}
\begin{eqnarray}
L^{SEP(j)}D_{n_{1}}(x_{1})D_{n_{2}}(x_{2}) 
& =  & n_{1}(2j-n_{2}) \left[ D_{n_{1}-1}(x_{1})D_{n_{2}+1}(x_{2}) - D_{n_{1}}(x_{1})D_{n_{2}}(x_{2})\right]  \\
& +  & (2j-n_{1})n_{2} \left[ D_{n_{1}+1}(x_{1})D_{n_{2}-1}(x_{2}) - D_{n_{1}}(x_{1})D_{n_{2}}(x_{2}) \right]  \nonumber
\end{eqnarray}
and the theorem is proved.
\cvd
\subsection{The Symmetric Inclusion Process, SIP($k$)} \label{sipprocess}
The Symmetric Inclusion Process with parameter $k >0$, denoted by SIP($k$), is a Markov jump process with unbounded state space where each site can have an arbitrary number of particles.
Again, we define the process on an undirected connected $G=(V,E)$ with $|V|=N$. Jumps occur at rate proportional to the number of particles in the departure and the arrival sites, as the generator describes:
\begin{equation}\label{sipg}
L^{SIP(k)}f(\textbf{x})=\sum_{\substack{1\le i < l \le N   \\  (i,l)\in E}} x_{i}(2k+x_{l})\left[ f(\textbf{x}^{i,l})-f(\textbf{x})\right]+x_{l}(2k+x_{i})\left[ f(\textbf{x}^{l,i})-f(\textbf{x})\right].
\end{equation}
Detailed balance is satisfied by a product measure with marginals given by identical Negative Binomial distributions with parameters $2k>0$ and $0 < p < 1$, i.e. with probability mass function
\be
\label{rho-neg-bin}
\rho(x) = {2k+x-1 \choose x} p^x (1-p)^{2k}\;,\qquad\qquad x\in\{0,1,\ldots\}\;.
\ee
The polynomials that are orthogonal with respect to the Negative Binomial distribution are the Meixner polynomials $M_{n}(x)$ with parameter $2k$, first introduced in \cite{Meixner}.
Choosing  in Eq. \eqref{hyperdiscrete}
\begin{equation}
\sigma(x) = x \qquad\qquad \tau(x) = 2k p -x(1-p) \qquad\qquad \lambda_n = n(1-p)
\end{equation}
we have that the Meixner polynomials are solution of the difference equation
\begin{equation}
\label{diff-eq-meix}
x\left[ M_{n}(x+1)-2M_{n}(x)+M_{n}(x-1)\right] + \left( 2kp-x+xp\right)  \left[ M_{n}(x+1)-M_{n}(x)\right] +n(1-p)M_{n}(x)=0\;.
\end{equation}
They satisfy the orthogonal relation 
\begin{equation}
\sum_{x=0}^{\infty}M_{n}(x)M_{m}(x)\rho(x)=\delta_{m,n}d_{n}^{2}
\end{equation}
with $\rho$ as in \eqref{rho-neg-bin} and norm in $\ell^2(\mathbb{N}_0, \rho)$ given by
\be
d_n^2 = \dfrac{n!\Gamma(2k+n)}{p^{n}\Gamma(2k)}
\ee
where $\Gamma(x)$ is the Gamma function.
As consequence of the orthogonality they satisfy the recurrence relation \eqref{expan2}
with 
\be
\alpha_n = \dfrac{p}{p-1} \qquad \qquad \beta_n = \dfrac{n+pn+2kp}{1-p}  \qquad \qquad \gamma_n = \dfrac{n(n-1+2k)}{p-1}
\ee
which then becomes
\be
\label{recu-meix}
xM_n(x) = \dfrac{p}{p-1}M_{n+1}(x) + \dfrac{n+pn+2kp}{1-p}M_n(x) + \dfrac{n(n-1+2k)}{p-1}M_{n-1}(x)\;.
\ee
Furthermore the raising operator in equation \eqref{raising-discrete} provides the identity 
%\be
%a f_n(x)  = -p(n+2k+x) f_n(x) + xf_n(x-1) 
%\ee
%with 
%\be
%c(n) =  -c
%\ee
%and thus we have the equation
\be
\label{raising-meix}
 [p(n+2k+x)] M_n(x) - xM_n(x-1)  =  pM_{n+1}(x) \;.
\ee

\subsubsection{Self-duality for SIP($k$)}
In analogy with the result for the Exclusion process it is possible to find a duality function for the Symmetric Inclusion Process in terms of the Meixner polynomials.
\begin{teor}
The SIP($k$) is a self-dual Markov process with self-duality function 
\begin{equation}
\label{meix}
D_{{\bf n}}({\bf x})= \prod_{i=1}^N \dfrac{\Gamma(2k)}{\Gamma(2k+n_i)}M_{n_i}(x_i)
\end{equation}
where $M_{n}(x)$ is the Meixner polynomial of degree $n$.
\end{teor}
\begin{oss}
The self-duality function can be rewritten in terms of hypergeometric function as
\begin{equation}
D_{\bf n}({\bf x})= \prod_{i=1}^N \mathstrut_2 F_1 \left( {\left. \genfrac{}{}{0pt}{} {-n_i,-x_i } { 2k }  \right\vert {1-\frac{1}{p}}} \right)\;.
\end{equation}
\end{oss}
\dim
As was done for the Exclusion Process we verify the self-duality relation in Equation \eqref{self-duall} for two sites, say $1$ and $2$. 
The action of the $SIP(k)$ generator working on the self-duality function for two sites is given by
\begin{eqnarray}
L^{SIP(k)}D_{n_{1}}(x_{1})D_{n_{2}}(x_{2}) 
& =  & x_{1}(2k+x_{2}) \left[ D_{n_{1}}(x_{1}-1)D_{n_{2}}(x_{2}+1) - D_{n_{1}}(x_{1})D_{n_{2}}(x_{2})\right]   \\
& +  & (2k+x_{1})x_{2} \left[ D_{n_{1}}(x_{1}+1)D_{n_{2}}(x_{2}-1) - D_{n_{1}}(x_{1})D_{n_{2}}(x_{2}) \right]\;. \nonumber
\end{eqnarray}
We rewrite this by factorizing site $1$ and $2$, i.e.
\begin{eqnarray}
\label{step2-meix}
L^{SIP(k)}D_{n_{1}}(x_{1})D_{n_{2}}(x_{2}) 
& =  & x_{1}D_{n_{1}}(x_{1}-1)(2k+x_{2})D_{n_{2}}(x_{2}+1)-x_{1}D_{n_{1}}(x_{1})(2k+x_{2})D_{n_{2}}(x_{2}) \\ 
& + &(2k+x_{1})D_{n_{1}}(x_{1}+1) x_{2}D_{n_{2}}(x_{2}-1)- (2k+x_{1}) D_{n_{1}}(x_{1})x_{2}D_{n_{2}}(x_{2}) \nonumber
\end{eqnarray} 
so that we now need an expression for the following terms:
\be
\label{3-termini-meix}
xD_{n}(x)\;, \qquad\quad xD_{n}(x-1)\;, \qquad\quad (2k+x)D_{n}(x+1)\;.
\ee
To get those, we first write the difference equation \eqref{diff-eq-meix}, the recurrence relation \eqref{recu-meix} and the raising operator equation \eqref{raising-meix} in terms of $D_{n}(x)$ using
\be
D_{n}(x)= \dfrac{\Gamma(2k)}{\Gamma(2k+n)}M_{n}(x)\;.
\ee   
Then the first term in \eqref{3-termini-meix} is simply obtained from the normalized recurrence relation, whereas the second and third terms are provided by simple algebraic manipulation of the normalized raising operator equation and the normalized difference equation.
We have, 
\begin{eqnarray}
xD_{n}(x)=\dfrac{p}{p-1}(2k + n)D_{n+1}(x)-\dfrac{n+p(n+2k)}{p-1}D_{n}(x)+ \dfrac{n}{p-1}D_{n-1}(x)\\
xD_{n}(x-1)=\dfrac{p}{p-1}(2k + n)D_{n+1}(x) -\dfrac{p}{p-1}(2k + 2n)D_{n}(x)+\dfrac{p}{p-1}nD_{n-1}(x)\\
(2k + x)D_{n}(x+1)=\dfrac{p}{p-1}(2k +n)D_{n+1}(x)-\dfrac{1}{p-1}(2k +2n)D_{n}(x)+\dfrac{1}{p-1}nD_{n-1}(x)\;.
\end{eqnarray}
These relations allow us to expand the SIP($k$) generator in Equation \eqref{step2-meix} as
\begin{align*}
L^{SIP}&D_{n_{1}}(x_{1})D_{n_{2}}(x_{2})  
\\ &=
\left[  \dfrac{p(2k + n_{1})}{p-1} D_{n_{1}+1}(x_1) - \dfrac{p(2k+2n_{1})}{p-1} D_{n_{1}}(x_1) + \dfrac{pn_{1}}{p-1} D_{n_{1}-1}(x_1) \right] \times 
\\ & 
\left[ \dfrac{p(2k +n_{2})}{p-1}D_{n_{2}+1}(x_{2})-\dfrac{2k +2n_{2}}{p-1}D_{n_{2}}(x_{2})+\dfrac{n_{2}}{p-1}D_{n_{2}-1}(x_{2}) \right] 
\\&
- \left[ \dfrac{p(2k + n_{1})}{p-1}D_{n_{1}+1}(x_{1})-\dfrac{n_{1}+p(n_{1}+2k)}{p-1}D_{n_{1}}(x_{1})+ \dfrac{n_{1}}{p-1}D_{n_{1}-1}(x_{1})        \right] \times
\\& 
\left[ \dfrac{p(2k + n_{2})}{p-1}D_{n_{2}+1}(x_{2})-\dfrac{n_{2}+p(n_{2}+2k)}{p-1}D_{n_{2}}(x_{2})+ \dfrac{n_{2}}{p-1}D_{n_{2}-1}(x_{2}) + 2k D_{n_{2}}(x_{2}) \right]     
\\&
+ \left[ \dfrac{p(2k + n_{2})}{p-1} D_{n_{2}+1}(x_2) - \dfrac{p(2k+2n_{2})}{p-1} D_{n_{2}}(x_2) + \dfrac{pn_{2}}{p-1} D_{n_{2}-1}(x_2) \right] \times 
\\ & 
\left[ \dfrac{p(2k +n_{1})}{p-1}D_{n_{1}+1}(x_{1})-\dfrac{2k +2n_{1}}{p-1}D_{n_{1}}(x_{1})+\dfrac{n_{1}}{p-1}D_{n_{1}-1}(x_{1}) \right] 
\\&
- \left[ \dfrac{p(2k + n_{2})}{p-1}D_{n_{2}+1}(x_{2})-\dfrac{n_{2}+p(n_{2}+2k)}{p-1}D_{n_{2}}(x_{2})+ \dfrac{n_{2}}{p-1}D_{n_{2}-1}(x_{2})        \right] \times
\\& 
\left[ \dfrac{p(2k + n_{1})}{p-1}D_{n_{1}+1}(x_{1})-\dfrac{n_{1}+p(n_{1}+2k)}{p-1}D_{n_{1}}(x_{1})+ \dfrac{n_{1}}{p-1}D_{n_{1}-1}(x_{1}) + 2k D_{n_{1}}(x_{1}) \right]\;.
\end{align*}

At this point it is sufficient to notice that the coefficients of products of polynomials with degree different than $n_1+n_2$ are all zero, so that we are left with  
\begin{eqnarray}
L^{SIP(k)}D_{n_{1}}(x_{1})D_{n_{2}}(x_{2}) 
& =  & n_{1}(2k+n_{2}) \left[ D_{n_{1}-1}(x_{1})D_{n_{2}+1}(x_{2}) - D_{n_{1}}(x_{1})D_{n_{2}}(x_{2})\right]  \\
& +  & (2k+n_{1})n_{2} \left[ D_{n_{1}+1}(x_{1})D_{n_{2}-1}(x_{2}) - D_{n_{1}}(x_{1})D_{n_{2}}(x_{2}) \right] \nonumber
\end{eqnarray}
and the theorem is proved.
\cvd

\subsection{The Independent Random Walker, IRW}
The Symmetric Independent Random Walkers, denoted IRW, is one of the simplest, yet non-trivial particle system studied in the literature. It consists of independent particles that perform a symmetric continuous time random walk at rate $1$.
The generator, defined on the undirected connected graph $G=(V,E)$ with $N$ vertices and edge set $E$, is given by
\begin{equation}\label{irwgen}
L^{IRW}f(\textbf{x})= \sum_{\substack{1\le i < l \le N   \\  (i,l)\in E}}  x_{i}\left[ f(\textbf{x}^{i,l})-f(\textbf{x})\right]+x_{l}\left[ f(\textbf{x}^{l,i})-f(\textbf{x})\right]\;.
\end{equation}
The reversible invariant measure is provided by a product of Poisson distributions with parameter $\lambda>0$,
 i.e. with probability mass function
\be
\label{rho-irw}
\rho(x) = \dfrac{e^{-\lambda} \lambda^{x}}{x!}\,, \qquad\qquad  x\in\mathbb{N}_0\;.
\ee
The orthogonal polynomials with respect the Poisson distribution are the Charlier polynomials $C_{n}(x)$ \cite{Cha} with parameter $\lambda$.
Choosing  in Eq. \eqref{hyperdiscrete}
\begin{equation}
\sigma(x) = x \qquad\qquad \tau(x) = \lambda-x \qquad\qquad \lambda_n = n
\end{equation}
we obtain the difference equation whose solution is $C_{n}(x)$
\begin{equation}\label{dec}
x\left[ C_{n}(x+1)-2C_{n}(x)+C_{n}(x-1)\right] + \left( \lambda -x \right)  \left[ C_{n}(x+1)-C_{n}(x)\right] +nC_{n}(x)=0\;.
\end{equation}
The orthogonal relation satisfied by Charlier polynomials is
\begin{equation}
\sum_{x=0}^{\infty}C_{n}(x)C_{m}(x)\rho(x)=\delta_{m,n}d_{n}^{2}
\end{equation}
where $\rho$ is given in \eqref{rho-irw} and the norm in $\ell^2(\mathbb{N}_0, \rho)$ is
\be
d_n^2 = n!\lambda^{-n}
\ee
As consequence of the orthogonality they satisfy the recurrence relation \eqref{expan2}
with 
\be
\alpha_n = -\lambda \qquad \qquad \beta_n = n+\lambda  \qquad \qquad \gamma_n = -n
\ee
which then becomes
\be \label{recc}
xC_n(x) = -\lambda C_{n+1}(x) + (n+\lambda)C_n(x) -nC_{n-1}(x)\;.
\ee
Furthermore, the raising operator in eq. \eqref{raising-discrete} provides
%\be
%a f_n(x)  = +\lambda f_n(x) - xf_n(x-1) 
%\ee
%with 
%\be
%c(n) =  -\lambda
%\ee
%and thus we have the equation
\be \label{raic}
\lambda C_n(x) - xC_n(x-1) =  \lambda C_{n+1}(x)\;.
\ee
\subsubsection{Self-duality of IRW}
As the following theorem shows, the self-duality relation is given by the Charlier polynomials themselves.
\begin{teor}
The IRW is a self-dual Markov process with self-duality function
\begin{equation}\label{cd}
D_{{\bf n}}({\bf{x}})=\prod_{i=1}^N C_{n_i}(x_i)
\end{equation}
where $C_{n}(x)$ is the Charlier polynomial of degree $n$.
\end{teor}
\begin{oss} Reading the Charlier polynomial as hypergeometric function, the duality function then becomes
\begin{equation}
D_{\textbf{n}}(\textbf{x})=\prod_{i=1}^N \mathstrut_2 F_0 \left( {\left. \genfrac{}{}{0pt}{} {-n_i,-x_i } { - }  \right\vert {-\frac{1}{\lambda}}} \right)\;.
\end{equation}
\end{oss}
\dim
It is clear from \eqref{cd} that the difference equations, the recurrence relations and the raising operator for $D_{n}(x)$ are respectively \eqref{dec}, \eqref{recc} and \eqref{raic}, that we rewrite as:
\begin{eqnarray}
D_{n}(x+1)=D_{n}(x)-\dfrac{n}{\lambda}D_{n-1}(x) \\
xD_{n}(x)=-\lambda D_{n+1}(x)+(n+\lambda)D_{n}(x)-nD_{n-1}(x)\\
xD_{n}(x-1)=\lambda D_{n}(x)-\lambda D_{n+1}(x).
\end{eqnarray} 
As done before, we use the two particles IRW generator in \eqref{irwgen} and the three equations above to check that the self-duality relation holds.
We have
\begin{align*}
L^{IRW}D_{n_{1}}(x_{1})D_{n_{2}}(x_{2}) 
& =  
x_{1}D_{n_{1}}(x_{1}-1)D_{n_{2}}(x_{2}+1)-x_{1}D_{n_{1}}(x_{1})D_{n_{2}}(x_{2}) \\
& + 
D_{n_{1}}(x_{1}+1)x_{2}D_{n_{2}}(x_{2}-1)-D_{n_{1}}(x_{1})x_{2}D_{n_{2}}(x_{2}) \\
& =
\left[ \lambda D_{n_{1}}(x_{1})-\lambda D_{n_{1}+1}(x_{1}) \right]\left[ D_{n_{2}}(x_{2}) - \dfrac{n_{2}}{\lambda}D_{n_{2}-1}(x_{2})\right]  \\ 
& - 
\left[-\lambda D_{n_{1}+1}(x_{1})+(n_{1}+ \lambda) D_{n_{1}}(x_{1}) -n_{1}D_{n_{1}-1}(x_{1}) \right] \left[ D_{n_{2}}(x_{2})\right] \\ 
& +
\left[D_{n_{1}}(x_{1}) - \dfrac{n_{1}}{\lambda}D_{n_{1}+1}(x_{2})\right] \left[ \lambda D_{n_{2}}(x_{2})-\lambda D_{n_{2}+1}(x_{2}) \right] \\ 
& -
\left[ D_{n_{1}}(x_{1}) \right] \left[-\lambda D_{n_{2}+1}(x_{2})+(n_{2}+ \lambda) D_{n_{2}}(x_{2}) -n_{2}D_{n_{2}-1}(x_{2}) \right]\;.
\end{align*}
After computing the products and suitable simplifications we get
\begin{eqnarray}
L^{IRW}D_{n_{1}}(x_{1}) D_{n_{2}}(x_{2})&=n_{1}[D_{n_{1}-1}(x_{1})D_{n_{2}+1}(x_{2})-D_{n_{1}}(x_{1})D_{n_{2}}(x_{2})]  \\ & \nonumber + n_2[ D_{n_{1}+1}(x_{1})D_{n_{2}-1}(x_{2})-D_{n_{1}}(x_{1})D_{n_{2}}(x_{2})]\;.
\end{eqnarray}
\cvd
\color{black}
%\textcolor{blue}{an alternative proof with Raising and Lowering operators.\\
%For Charlier polynomials the Raising and Lowering operators are respectively
%\begin{eqnarray}
%xC_{n}(x-1)-\lambda C_{n}(x)=-\lambda C_{n+1}(x) \\
%(xC_{n}(x-1)-xC_{n}(x)-n C_{n}(x))=(-\nabla -n)C_{n}(x)=nC_{n-1}(x)
%\end{eqnarray}
%Looking at the BEP generator, we see that we need and expression in terms of n's of
%\begin{equation}
%xC_{n}(x-1) \qquad xC_{n} \qquad C(x+1)
%\end{equation}
%We can use the Raising operator for the first one and the difference between lowering and raising for the middle one.\\
%How to find the third one?\\
%Facendo uno shift from $x$ to $x+1$ of the raising operator we get $xC_{n}(x)+C_{n}(x)-\lambda C_{n}(x+1)=-\lambda C_{n+1}(x+1)$ where expliciting for what we are looking for 
%\begin{equation*}
%\lambda C_{n}(x+1)=xC_{n}(x)+C_{n}(x)+\lambda C_{n+1}(x+1)
%\end{equation*}
%and we are stuck because we don't know how to write $C_{n+1}(x+1)$.\\
%On the other hand, if we use the difference equation we get
%\begin{equation*}
%\lambda C_{n}(x+1)=\lambda C_{n}(x)-nC_{n-1}(x)
%\end{equation*}
%Per questo mi sembra che per il BMP funzioni perchè con raising and lowering si riesce a ricavere la diff equation, mentre negli altri casi no.}
%
%
%
%

\section{Duality: proof of Theorem \ref{teo} part ii)}\label{duality}
In this last section we show two examples of duality: the initial process is an interacting diffusion, while the dual one is a jump process, which, in particular, turns out to be the SIP process introduced in section \ref{sipprocess}. We also show an example of duality for  a redistribution model of Kipnis-Marchioro-Presutti type.

\subsection{The Brownian Momentum Proces, BMP}

The Brownian Momentum Process (BMP) is a Markov diffusion process introduced in 
%\cite{BO} and 
\cite{GK}. 
%Think of the process described in Section \ref{bm} for $2$ sites now extended to $N$ sites and make $4k$ "vertical" copies of it. If sites $x_{i}$ and $ x_{j} $ are connected (in the sense that they share their energy and momentum) in the first horizontal layer, then sites $x_{i, \alpha}$ and $ x_{j, \beta} $ will be connected for every $\alpha $ and $\beta$.\\
On the undirected connected graph $G=(V,E)$ with $N$ vertices and edge set $E$,  the generator reads
%takes in account all the possible momentum shared among particles on the lattice graph:
\begin{equation}\label{genbm}
L^{BMP}f({\bf x})=\sum_{\substack{1\le i < l \le N   \\  (i,l)\in E}} 
\left( 
x_{i} \frac{\partial f}{\partial x_l}({\bf x})-   x_l \frac{\partial f}{\partial x_i}({\bf x})
\right)^{2} 
\end{equation}
where $f: \mathbb{R}^N \to \mathbb{R}$ is a function in the domain of the generator. A configuration is denoted by $\textbf{x}=(x_{i})_{i \in V}$ where $x_{i} \in \mathbb{R}$ has to be interpreted as a particle momentum. 
A peculiarity of this process regards its conservation law: if the process is started from the configuration ${\bf x}$ then $ ||{\bf x}||_2^2 = \sum_{i=1}^{N} x^{2}_{i} $ is constant during the evolution, i.e. the total kinetic energy is conserved. 
%Observe that the generator in \eqref{genbm} describes a rotation between the axis $x_{i}$ and $x_{j}$.
%In particular, it is the sum of all generators in \eqref{genbm} if sites $(i, \alpha)$ and $(j, \beta)$ are connected.\\

The stationary reversible measure of the BMP process is given by a family of product measures with marginals given by independent centered Gaussian random variables with variance $\sigma^2>0$. 
Without loss of generality we will identify a duality function related to Hermite polynomials in the case the variance equals $1/2$.
The case with a generic value of the variance would be treated in a similar way by using generalized Hermite polynomials.  
Choosing  in \eqref{hdf}
\begin{equation}
\sigma(x) = 1 \qquad\qquad \tau(x) = -2x \qquad\qquad \lambda_n = 2n
\end{equation}
we acquire the differential equation satisfied by the Hermite polynomials $H_{n}(x)$ \cite{NSU}
\begin{equation}\label{deh}
H^{''}_{n}(x)-2xH^{'}_{n}(x) +2nH_{n}(x)=0\;.
\end{equation}
The orthogonal relation they satisfy is
\begin{equation}
\int_{-\infty}^{+\infty}H_{n}(x)H_{m}(x)\rho(x)dx=\delta_{m,n}d_{n}^{2}
\end{equation}
with density function
\be
\rho(x) = \dfrac{e^{-x^{2}} }{\sqrt{\pi}}
\ee
and norm in $L^2(\mathbb{R}, \rho)$ given by
\be
d_n^2 = 2^{n}n!
\ee
As consequence of the orthogonality they satisfy the recurrence relation \eqref{expan2-cont}
with 
\be
\alpha_n = \frac{1}{2} \qquad \qquad \beta_n = 0  \qquad \qquad \gamma_n = n
\ee
which then can be written as
\be \label{rech}
H_{n+1}(x)-2xH_{n}(x)+2nH_{n-1}(x)=0\;.
\ee
Furthermore the raising operator in eq. \eqref{creation} provides 
%\be
%a f_n(x)  = \left( 2x-\partial_{x} \right) f_n(x) 
%\ee
%with 
%\be
%c(n) =  1
%\ee
%and thus we have the equation
\be \label{raih}
2xH_{n}(x)-H^{'}_n(x) = H_{n+1}(x).
\ee

\subsubsection{Duality between BMP and SIP($\frac14$)}
The following theorem has a similar version in \cite{CED} and it states the duality result involving the Hermite polynomials.
\begin{teor}
The BMP process is dual to the SIP($\frac14$) process through duality function
\begin{equation}
D_{\textbf{n}}(\textbf{x})=\prod_{i=1}^{N}\dfrac{1}{(2n_{i}-1)!!}H_{2n_{i}}(x_{i})
\end{equation}
where $H_{2n}(x)$ is the Hermite polynomial of degree $2n$.
\begin{oss} Reading the Hermite polynomial as hypergeometric function, the duality function then becomes
\begin{equation}
D_{\textbf{n}}(\textbf{x})=\prod_{i=1}^{N} \dfrac{1}{(2n_{i}-1)!!}(2x_{i})^{2n_{i}} \mathstrut_2 F_0 \left( {\left. \genfrac{}{}{0pt}{} {-n_{i},(-2n_{i}+1)/2 } { - }  \right\vert {-\frac{1}{x_{i}^{2}}}} \right)\;.
\end{equation}
\end{oss}
\end{teor}
\dim
%The orthogonality relation for the so defined duality function $D_{n}(x)$ is
%\begin{equation}
%\dfrac{1}{\sqrt{\pi}}\int_{- \infty}^{+\infty} D_{n}(x)D_{m}(x) e^{-x^2}=\dfrac{2^{n} n!}{(2n-1)!!} \delta_{n,m}
%\end{equation}
Although the proof in \cite{CED} can be easily adapted to our case, we show here an alternative proof that follows our general strategy of using the structural properties of Hermite polynomials. It is sufficient, as before, to show the duality relation in equation \eqref{du} for sites $1$ and $2$.
The action of the BMP generator on duality function reads
\begin{align}
L^{BMP}D_{n_{1}}(x_1)D_{n_{2}}(x_2)& =(x_{1}\partial_{x_{2}}-x_{2}\partial_{x_{1}})^2 D_{n_{1}}(x_1)D_{n_{2}}(x_2)
  \\& =
x_{1}^{2}D_{n_{1}}(x_1)D_{n_{2}}''(x_2) + D_{n_{1}}''(x_1)x_{2}^{2}D_{n_{2}}(x_2) -  x_{1}D'_{n_{1}}(x_1)D_{n_{2}}(x_2) \nonumber \\&
 -D_{n_{1}}(x_1)x_{2}D_{n_{2}}'(x_2) -2x_{1}D_{n_{1}}'(x_{1})x_{2}D_{n_{2}}'(x_{2}) \nonumber
\end{align}
where we use $\partial_{x_{i}}=\frac{\partial}{\partial x_{i}}$. %and the prime denotes derivative w.r.t the argument.
We now need the recurrence relation and the raising operator appropriately rewritten in term of the duality function in order to get suitable expression for
\begin{equation}
x^{2}D_{n}(x), \qquad\quad D_{n}''(x), \qquad\quad xD_{n}'(x).
\end{equation}
This can be done using 
\begin{equation}
D_{n}(x)=\dfrac{1}{(2n-1)!!}H_{2n}(x)
\end{equation}
so that
\begin{eqnarray}
x^{2}D_{n}(x)=\dfrac{1}{4}(2n+1)D_{n+1}(x)+\left( 2n+\dfrac{1}{2}\right)D_{n}(x)+2nD_{n-1}(x) \label{ddddd} \\
D^{''}_{n}(x)=8nD_{n-1}(x) \label{ddd} \\
xD^{'}_{n}(x)=2nD_{n}(x)+4nD_{n-1}(x)\label{dddd}
\end{eqnarray}
where \eqref{ddddd} is obtained from iterating twice the recurrence relation in \eqref{rech}, for equation \eqref{dddd} we combined \eqref{rech} and \eqref{raih} and then \eqref{ddd} is found from the differential equation \eqref{deh} using \eqref{dddd}.
Proceeding with the substitution into the generator we find
\begin{align*}
 L^{BMP}D_{n_{1}}(x_1)D_{n_{2}}(x_2) &=
 \left( \dfrac{1}{4}(2n_{1}+1)D_{n_{1}+1}(x_{1})+\left( 2n_{1}+\dfrac{1}{2}\right)D_{n_{1}}(x_{1})+2n_{1}D_{n_{1}-1}(x_{1}) \right) 8n_{2}D_{n_{2}-1}(x_{2}) 
  \nonumber \\&
  + 8n_{1}D_{n_{1}-1}(x_{1})\left( \dfrac{1}{4}(2n_{2}+1)D_{n_{2}+1}(x_{2})+\left( 2n_{2}+\dfrac{1}{2}\right)D_{n_{2}}(x_{2})+2n_{2}D_{n_{2}-1}(x_{2}) \right) \nonumber \\&
  - \left( 2n_{1}D_{n_{1}}(x_{1})+4n_{1}D_{n_{1}-1}(x_{1})\right) D_{n_{2}}(x_2) - D_{n_{1}}(x_1)\left( 2n_{1}D_{n_{2}}(x_{2})+4n_{2}D_{n_{2}-1}(x_{2})\right) \nonumber \\&
  - 2 \left( 2n_{1}D_{n_{1}}(x_{1})+4n_{1}D_{n_{1}-1}(x_{1})\right)\left( 2n_{2}D_{n_{2}}(x_{2})+4n_{2}D_{n_{2}-1}(x_{2})\right) \nonumber
\end{align*}
Finally, after appropriate simplification of the terms whose degree is different from $n_{1}+n_{2}$, we get
\begin{align}
L^{BMP}  D_{n_{1}}(x_1)D_{n_{2}}(x_2) & = 
(2n_{1}+1)2n_{2}\left[ D_{n_{1}+1}(x_1)D_{n_{2}-1}(x_2) - D_{n_{1}}(x_1)D_{n_{2}}(x_2)\right]  \\ &
+ 2n_{1}(2n_{2}+1)\left[ D_{n_{1}-1}(x_1)D_{n_{2}+1}(x_2) - D_{2n_{1}}(x_1)D_{n_{2}}(x_2)\right] \nonumber \\& =
L^{SIP}  D_{n_{1}}(x_1)D_{n_{2}}(x_2) \nonumber
\end{align}
which proves the theorem.

\cvd

\subsection{The Brownian Energy Process, BEP(k)}
We now introduce a process, known as Brownian Energy Process with parameter $k$, BEP($k$) in short notation, whose generator is
\begin{equation}\label{bepp}
L^{BEP(k)}f(\textbf{x})= \sum_{\substack{1\le i < l \le N   \\  (i,l)\in E}} \left[ x_{i}x_{j}\left( \dfrac{\partial}{\partial x_{i}} f(\textbf{x}) - \dfrac{\partial}{\partial x_{j}} f(\textbf{x})  \right)^{2}  + 2k(x_{i}-x_{j})\left( \dfrac{\partial}{\partial x_{i}} f(\textbf{x}) - \dfrac{\partial}{\partial x_{j}} f(\textbf{x})  \right)   \right].
\end{equation}
where $f: \mathbb{R}^N \to \mathbb{R}$ is in the domain of the generator and  $\textbf{x}=(x_{i})_{i \in V}$ denotes a configuration of the process with $x_{i} \in \mathbb{R^{+}}$ interpreted as a particle energy.
The generator in \eqref{bepp} describes the evolution of particle system that exchange their (kinetic) energies. It is easy to verify that the total energy of the system $ \sum_{i=1}^{N} x_{i} $ is conserved by the dynamic.\\
In \cite{GKRV} it was shown that the BEP($k$) can be obtained from the BMP process once $4k \in \mathbb{N}$ vertical copies of the graph $G$ are introduced. Under these circumstances denoting $z_{i, \alpha}$ the momentum of the $i^{th}$ particle at the $\alpha^{th}$ level, the kinetic energy per (vertical) site is
\begin{equation}
x_{i}= \sum_{\alpha=1}^{4k} z_{i, \alpha}^{2} \;.
\end{equation}
If we use the above change of variable in the generator of such BMP process on the ladder graph with $4k$ layers, the generator of the BEP($k$) is revealed.
%\begin{figure}[hbtp] 
%\centering
%\includegraphics[scale=0.4]{../../../Dropbox/Tesi/Dottorato/Lattice.jpg}
%\caption{Lattice}
%\label{fig:lat}
%\end{figure}

The stationary measure of the BEP($k$) process is given by a product of independent Gamma distribution with shape parameter $2k$ and scale parameter $\theta$, i.e. with Lebesgue probability mass function
\be\label{BEPmassfct}
\rho(x) = \dfrac{x^{2k-1}e^{-\frac{x}{\theta}}} {\Gamma(2k) \theta^{k}}.
\ee 
Without loss of generality we can set $\theta=1$ so that the polynomials orthogonal with respect to the Gamma distribution are the generalized Laguerre polynomials $L^{(2k-1)}_{n}(x)$ \cite{NSU}. \\
Choosing 
\begin{equation}
\sigma(x) = x \qquad\qquad \tau(x) = 2k-x \qquad\qquad \lambda_n = n
\end{equation}
the differential equation whose generalized Laguerre polynomials are solution becomes
\begin{equation}\label{del}
x\frac{d^{2}}{dx^{2}}L^{(2k-1)}_{n}(x)+(2k-x)\frac{d}{dx}L^{(2k-1)}_{n}(x) +nL^{(2k-1)}_{n}(x)=0
\end{equation}
The orthogonal relation they satisfy is 
\begin{equation}
\int_{0}^{+\infty}L^{(2k-1)}_{n}(x)L^{(2k-1)}_{m}(x)\rho(x)dx=\delta_{m,n}d_{n}^{2}
\end{equation}
with mass function as in \eqref{BEPmassfct} with $\theta=1$ and norm in $L^2(\mathbb{R^{+}}, \rho)$ given by
\be
d_n^2 = \dfrac{\Gamma(n+2k)}{n! \Gamma(2k)}\;.
\ee
As consequence of the orthogonality, they satisfy the recurrence relation \eqref{expan2-cont}
with 
\be
\alpha_n = -(n+1) \qquad \qquad \beta_n = 2n+2k  \qquad \qquad \gamma_n = -(n+2k-1)
\ee
which then can be written as
\be \label{recl}
xL_{n}^{(2k-1)}(x)=-(n+1)L_{n+1}^{(2k-1)}(x)+(2n+2k)L_{n}^{(2k-1)}(x)-(n+2k-1)L_{n-1}^{(2k-1)}(x)
\ee
Furthermore, the raising operator in eq. \eqref{creation} is given by 
%\be
%a f_n(x)  = \left( x\partial_{x} +2k-x+n \right) f_n(x) 
%\ee
%with 
%\be
%c(n) =  n+1
%\ee
%and thus we have the equation
\be \label{rail}
(2k-x+n)L_{n}^{(2k-1)}(x)+x\frac{d}{dx}L_{n}^{(2k-1)}(x)=(n+1)L_{n+1}^{(2k-1)}(x)\;.%=nL_{n}^{2k-1}(z)-(n+2k-1)L_{n-1}^{2k-1}(z)
\ee
%where we used \eqref{recl} for the last equality above.

\subsubsection{Duality between BEP($k$) and SIP($k$)}
The duality relation for the Brownian energy process with parameter $k$ is stated below.
\begin{teor}
The BEP($k$) process and the SIP($k$) process are dual via 
\begin{equation}\label{bepdu}
D_{\textbf{n}}(\textbf{x})=\prod_{i=1}^{N} \dfrac{n_{i}! \, \Gamma(2k)}{\Gamma(2k+n_{i})}L_{n_{i}}^{\left( 2k-1 \right) }(x_{i})
\end{equation}
where $L_{n}^{\left( 2k-1 \right)}(x)$ is the generalized Laguerre polynomial of degree $n$.
\end{teor}
\begin{oss} The factor $\Gamma(2k)$ in \eqref{bepdu} is not crucial to assess a duality relation, but it allows to write the duality function as the hypergeometric function
\begin{equation}
D_{\textbf{n}}(\textbf{x})= \prod_{i=1}^{N} \mathstrut_1 F_1 \left( {\left. \genfrac{}{}{0pt}{} {-n_{i}} { 2k}  \right\vert {x_{i}}} \right)\;.
\end{equation}
\end{oss}
\dim %Note that the last equality only holds for duality purpose.\\
%The orthogonality relation for the so defined duality function $D_{n}(x)$ is
%\begin{equation}
%\int_{0}^{+\infty} D_{n}(x)D_{m}(x) \dfrac{x^{\alpha} e^{-x}}{\Gamma(\alpha +1)}=\dfrac{n! \alpha!}{(n+\alpha)!} \delta_{n,m}
%\end{equation}
As in the previous cases we notice that the proof can be shown for sites $1$ and $2$ only, in which case the generator of the BEP acts on 
\begin{align}
L^{BEP(k)} D_{n_{1}}(x_{1})D_{n_{2}}(x_{2})=  & \left[ x_{1}x_{2}\left( \partial_{x_{1}} - \partial_{x_{2}} \right)^{2}- 2k(x_{1}-x_{1})(\partial_{x_{1}}-\partial_{x_{2}}) \right] D_{n_{1}}(x_{1})D_{n_{2}}(x_{2})  \nonumber \\ =& 
(x_{1}\partial_{x_{1}}^{2}+2k\partial_{x_{1}})D_{n_{1}}(x_{1})x_{2}D_{n_{2}}(x_{2})+x_{1}D_{n_{1}}(x_{1})(x_{2}\partial_{x_{2}}^{2}+2k\partial_{x_{2}})D_{n_{2}}(x_{2})  \nonumber \\ - & 
x_{1}\partial_{x_{1}}D_{n_{1}(x_{1})}(x_{2}\partial_{x_{2}}+2k)D_{n_{2}}(x_{2})-(x_{1}\partial_{x_{1}}+2k)D_{n_{1}}(x_{1})x_{2}\partial_{x_{2}}D_{n_{2}}(x_{2}) \nonumber
\end{align}
We seek an expression for
\begin{equation}
x \partial_{x}^{2}D_{n}+2k \partial_{x}D_{n}, \qquad xD_{n}, \qquad x \partial_{x}D_{n}
\end{equation}
that can easily be obtained rewriting \eqref{del}, \eqref{recl} and \eqref{rail} for the duality function, using 
\begin{equation}
D_{n}(x)=\dfrac{n! \, \Gamma(2k)}{\Gamma(2k+n)}L_{n}^{\left( 2k-1 \right) }(x)
\end{equation}
so that, after simple manipulation
\begin{eqnarray} 
\label{questaqui}
xD^{''}_{n}(x)+(2k -x)D^{'}_{n}(x)+nD_{n}(x)=0 \\
xD_{n}(x)=-(n+2k)D_{n+1}(x)+(2n+2k)D_{n}(x)-nD_{n-1}(x)\\
xD_{n}^{'}(x)=nD_{n}(x)-nD_{n-1}(x)\;. \label{raiL}
\end{eqnarray}
Note that plugging \eqref{raiL} into the difference equation \eqref{questaqui}, we get 
$xD_{n}^{''}(x)+2kD^{'}_{n}(x)=-nD_{n-1}(x)$. \\
Let's now use these information to write explicitly the BEP($k$) generator. 
\begin{align*}
L^{BEP(k)} D_{n_{1}}(x_{1})D_{n_{2}}(x_{2}) & =
\left[ -n_{1}D_{n_{1}-1}(x_{1})\right] \left[ -(2k+n_{2})D_{n_{2}+1}(x_{2})+(2n_{2}+2k)D_{n_{2}}(x_{2})-n_{2}D_{n_{2}-1}(x_{2})\right] \\& + \left[ -(2k+n_{1})D_{n_{1}+1}(x_{1})+(2n_{1}+2k)D_{n_{1}}(x_{1})-  n_{1}D_{n_{1}-1}(x_{1}) \right]  \left[ -n_{2}D_{n_{2}-1}(x_{2}) \right]  \\& -
 \left[ n_{1}D_{n_{1}}(x_{1})-n_{1}D_{n_{1}-1}(x_{1})\right] \left[(n_{2}+2k)D_{n_{2}}(x_{2})-n_{2}D_{n_{2}-1}(x_{2}) \right]  \\& - \left[(n_{1}+2k)D_{n_{1}}(x_{1})-n_{1}D_{n_{1}-1}(x_{1}) \right]\left[ n_{2}D_{n_{2}}(x_{2})-n_{2}D_{n_{2}-1}(x_{2})\right] \;.
\end{align*}
Expanding products in the above expression we find
\begin{align*}
 L^{BEP(k)} D_{n_{1}}(x_{1})D_{n_{2}}(x_{2}) & = n_{1}D_{n_{1}-1}(x_{1})(n_{2}+2k)D_{n_{2}+1}(x_{2})+(n_{1}+2k)D_{n_{1}+1}(x_{1})n_{2}D_{n_{2}-1}(x_{2}) \\& +n_{1}D_{n_{1}}(x_{1})(n_{2}+2k)D_{n_{2}}(x_{2})+(n_{1}+2k)D_{n_{1}+1}(x_{1})n_{2}D_{n_{2}-1}(x_{2})\\& 
  + D_{n_{1}-1}(x_{1})D_{n_{2}}(x_{2})\left[ -n_{1}(2n_{2}+2k) +n_{1}(n_{2}+2k) + n_{1}n_{2}\right]  \\& + D_{n_{1}}(x_{1})D_{n_{2}-1}(x_{2})\left[ -(2n_{1}+2k)n_{1} +(n_{1}+2k)n_{2} +n_{1}n_{2}\right] \\&  
  + D_{n_{1}-1}(x_{1})D_{n_{2}-1}(x_{2})\left[n_{1}n_{2}+n_{1}n_{2}-n_{1}n_{2}-n_{1}n_{2}\right] \;.
\end{align*}
Noticing that the coefficients of the last three lines are zeros, we finally get 

\begin{align}
\label{plutone}
L^{BEP(k)} D_{n_{1}}(x_{1})D_{n_{2}}(x_{2})& =n_{1}(n_{2}+2k)\left[ D_{n_{1}-1}(x_{1})D_{n_{2}+1}(x_{2})  -D_{n_{1}}(x_{1})D_{n_{2}-1}(x_{2})\right] \\&
 + (n_{1}+2k)n_{2} \left[ D_{n_{1}+1}(x_{1})D_{n_{2}-2}(x_{2})-D_{n_{1}}(x_{1})D_{n_{2}}(x_{2}) \right] \nonumber \\ & = 
 L^{SIP(k)}D_{n_{1}}(x_{1})D_{n_{2}}(x_{2}) \nonumber
\end{align}
where $L^{SIP(k)}$ works on the dual variables $(n_{1},n_{2})$.
\cvd

\subsection{The Kipnis-Marchioro-Presutti process, KMP(k)}
The KMP model was first introduced by Kipnis, Marchioro and Presutti  \cite{KMP} in 1982 as a model of heat conduction that was solved by using a dual process.
It is a stochastic model where a continuous non-negative variable (interpreted as energy) is uniformly redistributed among two random particles on a lattice at Poisson random times.
A general version with parameter $k$, that we shall call KMP($k$) was defined in \cite{CGGR2},
by considering a redistribution rule where a fraction $p$ of the total energy is assigned to one particle and the remaining fraction ($1-p$) to the other particle, with $p$ a Beta($2k,2k$) distributed random variable. Thus the case $k=1/2$ corresponds to the original KMP model. 
In \cite{CGGR2} it was shown that KMP($k$) is in turn related to the Brownian Energy Process with parameter $k$, as it can be obtained from the BEP($k$)
via a procedure called ``instantaneous thermalization''. If the spatial setting remains as described in the previous sections, the generator of the KMP($k$) process is
\begin{equation}
L^{KMP(k)}f({\bf x})=\sum_{\substack{1\le i < l \le N   \\  (i,l)\in E}} 
\int_{0}^{1} \left[  f\left( x_{1},\ldots ,x_{i-1},p(x_{i}+x_{i+1}),(1-p)(x_{i}+x_{i+1}),x_{i+2},\ldots x_{N} \right) -f\left( \bf x\right) \right]  \nu_{2k}(p) dp
\end{equation}
where $\nu_{2k}(p)$ is the  density function of the Beta distribution with parameters ($2k$, $2k$), i.e.
\begin{equation}
\nu_{2k}(p)=\dfrac{p^{2k-1}(1-p)^{2k-1} \Gamma(4k)}{\Gamma(2k)\Gamma(2k)} \;, \qquad p \in (0,1)\;.
\end{equation}

The dual process of KMP($k$) \cite{CGGR2} is generated by 
\begin{equation}
L^{\text{dual-}KMP(k)}f({\bf n})=\sum_{\substack{1\le i < l \le N   \\  (i,l)\in E}} 
\sum_{r=0}^{n_{i}+n_{i+1}} \left[  f\left( n_{1},\ldots ,n_{i-1},r,n_{i}+n_{i+1}-r,n_{i+2},\ldots x_{N} \right) -f\left( \bf n\right) \right]  \mu_{2k}(r\mid n_{i}+n_{i+1})
\end{equation}
where $\mu_{2k}(r|C)$ is the mass density function of the Beta Binomial distribution with parameter ($C$, $2k$, $2k$), i.e.
\begin{equation}
\mu_{2k}(r \mid C)=\dfrac{ \binom {2k+r-1} {r} \binom {2k+C-r-1} {C-r}}{\binom {4k+C-1} {C}}\;, \qquad r\in \left\lbrace 0,1,\ldots, C \right\rbrace  \;.
\end{equation}
This generator is the result of a thermalized limit of the SIP($k$) \cite{CGGR2}.
Our last theorem is stated below.
\begin{teor}
The KMP($k$) process duality relation with its dual is established via duality function 
\begin{equation}\label{kmpdu}
D_{\textbf{n}}(\textbf{x})=\prod_{i=1}^{N} \dfrac{n_{i}! \, \Gamma(2k)}{\Gamma(2k+n_{i})}L_{n_{i}}^{\left( 2k-1 \right) }(x_{i})
\end{equation}
where $L_{n}^{\left( 2k-1 \right)}(x)$ is the generalized Laguerre polynomial of degree $n$.
\end{teor}
\dim
As expected, the duality function is the same as the one for the BEP($k$) and SIP($k$) duality relation. This shouldn't surprise since BEP($k$) and SIP($k$) are dual through duality function \eqref{kmpdu} and the thermalization limit doesn't affect the duality property. Indeed, considering two graph vertices, one has from \cite{CGGR2}
\be
L^{KMP(k)}f({x_1,x_2}) = \lim_{t\to\infty} (e^{tL^{BEP(k)}} -  I )f(x_1,x_2) 
\ee
and
\be
L^{\text{dual-}KMP(k)}f({n_1,n_2}) = \lim_{t\to\infty} (e^{tL^{SIP(k)}} - I)f(n_1,n_2)\;. 
\ee
Thus, combining the previous two equations and \eqref{plutone}, the claim follows.
\cvd

\noindent
\textbf{Acknowledgments.} This research was supported by the Italian Research Funding Agency (MIUR) through FIRB project
grant n. RBFR10N90W and in part by the National Science Foundation under Grant No. NSF PHY11-25915.
We acknowledge a useful discussion on the topic of this paper with C\'edric Bernardin during the trimester “Disordered systems, random spatial processes and their applications” that was held at the Institute Henri Poincar\'e.


\begin{thebibliography}{empty}


\bibitem{BarCor16} \label{BarCor16} G. Barraquand, I. Corwin. The $ q $-Hahn asymmetric exclusion process. The Annals of Applied Probability 26.4, 2304--2356 (2016).

\bibitem{BelSch15} \label{BelSch15} V. Belitsky, G.M. Schütz. Self-duality for the two-component asymmetric simple exclusion process. Journal of Mathematical Physics 56.8, 083302 (2015).

\bibitem{BelSch16} \label{BelSch16} V. Belitsky, G.M. Schütz. Self-duality and shock dynamics in the $ n $-component priority ASEP. Preprint arXiv:1606.04587 (2016).

\bibitem{CED} \label{CED} C. Bernardin. Superdiffusivity of asymmetric energy model in dimensions $1$ and $2$. Journal of Mathematical Physics 49.10 103301 (2008).

\bibitem{BO} \label{BO} C. Bernardin, S. Olla. Fourier’s law for a microscopic model of heat conduction. Journal of Statistical Physics 121.3--4, 271--289 (2005).

\bibitem{BorCor13} \label{BorCor13} A. Borodin, I. Corwin. Discrete time q-TASEPs. International Mathematics Research Notice Issue 2,  499--537 (2013).

\bibitem{BorCorGor16} \label{BorCorGor16} A. Borodin, I. Corwin, V. Gorin. Stochastic six-vertex model. Duke Mathematical Journal 165.3, 563--624 (2016).

\bibitem{BCPS15} \label{BCPS15} A. Borodin, I. Corwin, L. Petrov, T. Sasamoto. Spectral theory for the-Boson particle system. Compositio Mathematica 151.01, 1--67 (2015).

\bibitem{BorCorSa14} \label{BorCorSa14} A. Borodin, I. Corwin, T. Sasamoto. From duality to determinants for q-TASEP and ASEP. The Annals of Probability 42.6, 2314--2382 (2014).


\bibitem{CGGR2} \label{CGGR2} G. Carinci, C. Giardinà, C. Giberti, F. Redig. Duality for stochastic models of transport. Journal of Statistical Physics 152.4, 657--697 (2013).

\bibitem{CGGR} \label{CGGR} G. Carinci, C. Giardinà, C. Giberti, F. Redig. Dualities in population genetics: a fresh look with new dualities. Stochastic Processes and their Applications 125.3, 941--969 (2015).

\bibitem{CGRT} \label{CGRT} G. Carinci, C. Giardinà, F. Redig, T. Sasamoto.
A generalized Asymmetric Exclusion Process with $U_q(\mathfrak{sl}_2)$ stochastic duality. Probability Theory and Related Fields, 1--47 (2014).

\bibitem{CGRS16} \label{CGRS16} G. Carinci, C. Giardinà, F. Redig, T. Sasamoto. Asymmetric Stochastic Transport Models with ${\mathscr{U}}_q (\mathfrak{su}(1, 1))$ Symmetry. Journal of Statistical Physics 163.2, 239--279 (2016).

\bibitem{Cha} \label{Cha}
C.V.L. Charlier. \"Uber die darstellung willk\"urlicher Functionen. Arkiv f\"or Matematik, Astronomi och Fysik,  2.20, 1--35 (1905--1906).

\bibitem{TSC} \label{TSC} T.S. Chihara. \emph{An Introduction to orthogonal polynomials}. Gordon and Breach (1978).

%\bibitem{Fed}
%E. Fedheim: On a system of orthogonal polynomials associated with a distribution of Stieltjes type, Acad. Sci. URSS, 31 528-533 (1941)


\bibitem{Cor15} \label{Cor15} I. Corwin. The q-Hahn Boson process and q-Hahn TASEP. International Mathematics Research Notices Issue 14, 5577--5603 (2014).

\bibitem{CorPet16} \label{CorPet16} I. Corwin, L. Petrov. Stochastic higher spin vertex models on the line. Communications in Mathematical Physics 343.2, 651--700 (2016).

\bibitem{CST16} \label{CST16} I. Corwin, H. Shen, L-C Tsai. ASEP(q,j) converges to the KPZ equation. Preprint arXiv:1602.01908 (2016).

\bibitem{DP06}	\label{DP06}
A. De Masi, E. Presutti. {\em Mathematical methods for hydrodynamic limits}. Springer (2006).

%\bibitem{CHI} \label{CHI} C. Franceschini. Duality theory for stochastic processes with multiple conservation laws. (2014).

\bibitem{GK} \label{GK} C. Giardinà, J. Kurchan. The Fourier law in a momentum-conserving chain. Journal of Statistical Mechanics: Theory and Experiment 2005.05,  P05009 (2005).


\bibitem{GKR} \label{GKR} C. Giardinà, J. Kurchan, F. Redig. Duality and exact correlations for a model of heat conduction. Journal of Mathematical Physics 48.3, 033301 (2007).

\bibitem{GKRV} \label{GKRV} C. Giardinà, J. Kurchan, F. Redig, K. Vafayi. Duality and hidden symmetries in interacting particle systems. Journal of Statistical Physics 135.1, 25--55  (2009).
 

\bibitem{GRV10} \label{GRV10} C. Giardinà, F. Redig, K. Vafayi. Correlation inequalities for interacting particle systems with duality. Journal of Statistical Physics 141.2, 242--263 (2010).

\bibitem{GwaSpo92} \label{GwaSpo92} L. Gwa, H. Spohn. Bethe solution for the dynamical-scaling exponent of the noisy Burgers equation. Physical Review A 46.2, 844 (1992).

\bibitem{ImSa11} \label{ImSa11} T. Imamura, T. Sasamoto. Current moments of 1D ASEP by duality. Journal of Statistical Physics 142.5, 919--930 (2011).

\bibitem{JK} \label{JK} S. Jansen, N. Kurt. On the notion(s) of duality for Markov processes, Probability Surveys 11, 59-120 (2014).

\bibitem{KMP}\label{KMP} C. Kipnis, C. Marchioro, E. Presutti. Heat flow in an exactly solvable model. Journal of Statistical Physics 27.1, 65--74 (1982).

\bibitem{KLS} \label{KLS} R. Koekoek, P.A. Lesky, R.F. Swarttouw.  \emph{ Hypergeometric Orthogonal Polynomials and their $q-$Analogues}, Springer (2010).

\bibitem{Koornwinder} \label{Koornwinder} T.H. Koornwinder, Lowering and raising operators for some special orthogonal polynomials. Preprint arXiv math/0505378 (2005).

\bibitem{Kra}\label{Kra}
M. Krawtchouk. Sur une généralisation des polynomes d'Hermite. Comptes Rendus 189.620--622 (1929).

\bibitem{Kuan15}\label{Kuan15}
J. Kuan. Stochastic duality of ASEP with two particle types via symmetry of quantum groups of rank two. Journal of Physics A: Mathematical and Theoretical 49.11, 115002 (2016).

\bibitem{Kuan16}\label{Kuan16}
J. Kuan. A Multi-species ASEP (q, j) and q-TAZRP with Stochastic Duality. Preprint arXiv:1605.00691 (2016).

\bibitem{Kuan17}\label{Kuan17}
J. Kuan. An algebraic construction of duality functions for the stochastic $ U_q (A_n^{(1)})$ vertex model and its degenerations. Preprint arXiv:1701.04468 (2017).

\bibitem{Liggett} \label{Liggett} T. M. Liggett. \emph{Interacting particles systems} Springer (1985).

\bibitem{Meixner}\label{Meixner}
J. Meixner. Orthogonale Polynomsysteme mit einer besonderen Gestalt der erzeugenden Funktion. Journal of the London Mathematical Sociaty 1.1, 6--13 (1934).

\bibitem{M}\label{M}M. M\"ohle. The concept of duality and applications to Markov processes arising in neutral population genetics models. Bernoulli 5.5, 761-–777 (1999).

\bibitem{NSU}\label{NSU} A.F. Nikiforov, S.K. Suslov, V.B. Uvarov. \emph{Classical Orthogonal
Polynomials of a Discrete Variable}, Springer-Verlag, Berlin (1991).


\bibitem{Ohku16} \label{Ohku16} J. Ohkubo. On dualities for SSEP and ASEP with open boundary conditions. Journal of Physics A: Mathematical and Theoretical, (2017).

%\bibitem{STS} \label{STS} T. Sasamoto, H. Spohn. One-dimensional Kardar-Parisi-Zhang equation: an exact solution and its universality. Physical review letters 104.23, 230602 (2010).

\bibitem{W} \label{W} W. Schoutens. \emph{Stochastic Processes and Orthogonal Polynomials}, Springer (2000).

\bibitem{Sch97} \label{Sch97} G.M. Schütz. Duality relations for asymmetric exclusion processes. Journal of Statistical Physics 86.5, 1265--1287 (1997).

\bibitem{Schutz-Sandow94} \label{Schutz-Sandow94} G.M. Schütz, S. Sandow. Non-Abelian symmetries of stochastic processeHeat flow in an exactly solvable mods: Derivation of correlation functions for random-vertex models and disordered-interacting-particle systems. Physical Review E 49.4, 2726 (1994).

\bibitem{Spi70} \label{Spi70} F. Spitzer. Interaction of Markov processes. 
Advances in Mathematics 5.2, 246--290 (1970).

\bibitem{Spohn} \label{Spohn} H. Spohn. Long range correlations for stochastic lattice gases in a non-equilibrium steady state. Journal of Physics A: Mathematical and General 16.18, 4275 (1983).






\end{thebibliography}
\end{document}